\documentclass[11pt]{report}

\usepackage[ngerman]{babel}
\usepackage{latexsym} 
\usepackage{amsfonts} 
\usepackage{amsmath}
\usepackage{pstricks}
\usepackage{graphics}

\newcommand{\bq}{\begin{eqnarray}}
\newcommand{\eq}{\end{eqnarray}}
\newcommand{\bew}{{\bf Beweis: }}       
\newcommand{\qed}{\hfill \mbox{$\Box$}}

\begin{document}
\title{\huge{Verteilung des Geschlechts zuf"allig gew"ahlter Origamis}\\[2cm]}
\author{{\LARGE{Diplomarbeit}} \\ \\ {\large vorgelegt von:}\\ \Large{Sabine Lechner}\\[2cm]}
\date{{\large angefertigt} \\ \Large{ am Mathematischen Institut der \\ Albert-Ludwig-Universit"at Freiburg \\ Juni 2008 \\[1.5cm]
Betreuer: apl. Prof. Dr. Jan-Christoph Schlage-Puchta}}
\maketitle
\newpage

\setlength{\parskip}{0.12cm}
\pagenumbering{roman}
\tableofcontents



\newpage
\setlength{\parskip}{0.7cm}
\pagenumbering{arabic}
\chapter*{Einleitung}
\addcontentsline{toc}{chapter}{Einleitung}
Es ist ein wichtiges Forschungsziel in der algebraischen Geometrie und der komplexen Differentialgeometrie, die Geometrie des Modulraums $M_g$ der Kurven vom Geschlecht $g$ zu verstehen. Ein Ansatz ist es, algebraische Kurven dieses Modulraums, die das Bild einer Teichm"ullerkreisscheibe sind, zu untersuchen. Eine spezielle Klasse dieser Teichm"ullerkurven stellen die Origamikurven dar. Ein Origami (oft auch "'square tiled surface"' genannt) ist eine topologische Fl"ache vom Geschlecht $g$, die entsteht, indem man endlich viele Einheitsquadrate entlang der Kanten so verklebt, dass die resultierende Fl"ache kompakt und zusammenh"angend ist. Durch die zu solch einer Fl"ache geh"orende Translationsstruktur - eine komplexe Struktur, deren "Ubergangsfunktionen Translationen sind - erh"alt man eine Riemannsche Fl"ache vom Geschlecht $g$, also einen Punkt im Modulraum $M_g$. Eine spezielle Form von Origamis wurde schon von W. Thurston in seiner Arbeit "uber die Diffeomorphismen von Fl"achen in \cite{thur} definiert, die dann sp"ater von W. Veech in \cite{veech} wieder aufgenommen wurde. 

In der vorliegenden Arbeit wird das Ziel verfolgt, die Verteilung des Geschlechts $g$, ermittelt mit Hilfe der Eulerschen Formel $\chi=E-K+F$, der durch ein Origami enstehenden Fl"ache $X$, zu berechnen. Die Idee ist hierbei, die Anzahl der Ecken $E$ von $X$ "uber die Anzahl der Bahnen der Operation
\begin{eqnarray*} 
\phi: \left\langle \sigma\tau\right\rangle \times ([4]\times[n]) \rightarrow ([4]\times[n])\nonumber \\
\phi((\sigma\tau)^k,(i,j))=(\tilde{\pi}^k(i),(f_k)_{\tilde{\pi}^k(i)}(j))
\end{eqnarray*}
mit $\sigma$ und $\tau \in C_4\wr S_n$ zu berechnen. Hierf"ur werden in Kapitel 1 zun"achst die grundlegenden Begriffe zu Origamis und zum Kranzproduktes (bezeichnet mit $\wr$) eingef"uhrt. Die Wahrscheinlichkeit  wird mit Hilfe der Formel 
\[P(\# \mbox{Bahnen} =k)=\sum_{\pi:\pi' \ hat \ k \ Zykel} P(\sigma\tau=\pi)=\frac{1}{|C_4\wr S_n|}\sum_D \frac{\chi^D(\sigma)\chi^D(\tau)\chi^D(\pi)}{\chi^D(1)}\]
berechnet, wobei $\pi=(\tilde{\pi};\pi_1,\pi_2,\pi_3,\pi_4) \in C_u\wr S_n$, $\pi'=\pi_1\circ\pi_2\circ\pi_3\circ\pi_4$ und die Summe "uber alle irreduziblen Darstellungen  D von $C_4\wr S_n$ l"auft. Die daf"ur n"otigen Kenntnisse der Darstellungstheorie werden in Kapitel 2 behandelt. In Kapitel 3 wird die Idee und die vorangegangene Theorie zusammengef"uhrt, um die Verteilung des Geschlechts zu berechnen. Kapitel 4 soll einen kurzen Einblick geben, womit sich die Forschung rund um Origamis gerade besch"aftigt.

An dieser Stelle m"ochte ich mich bei all denen bedanken, die mir bei der Abfassung dieser Arbeit geholfen haben.
Mein besonderer Dank gilt dabei Herrn apl. Prof. Dr. Jan-Christoph Schlage-Puchta, der mich an ein reizvolles Teilgebiet der Mathematik herangef"uhrt hat. Zu gro"sem Dank bin ich auch der Arbeitsgruppe Origamis in Karlsruhe verpflichtet, die mir die M"oglichkeit gab, "uber mein Diplomarbeitsthema auf ihrem Weihnachtsworkshop zur (algebraischen) Geometrie und Zahlentheorie vorzutragen.

\chapter{Einf"uhrung}
\section{Origamis}

Origamis lassen sich auf verschiedene Weisen beschreiben, am anschaulichsten ist jedoch die folgende Definition.

{\bf Definition 1.1.1}\newline
Ein {\it Origami} erh"alt man durch die Verklebung endlich vieler Kopien des euklidischen Einheitsquadrates. Diese Verklebungen richten sich nach den folgenden Regeln:
\begin{list}{$\bullet$}{\setlength{\topsep}{0cm} \setlength{\parskip}{0.4cm}}
\item Jede rechte Kante wird mit einer linken Kante verklebt.
\item Jede obere Kante wird mit einer unteren Kante verklebt.
\item Die so entstehende abgeschlossene Fl"ache $X$ ist zusammenh"angend. 
\end{list}

Die Namensgebung Origami geht zur"uck auf Pierre Lochak. In seiner Arbeit \cite{Lochak} verwendet er erstmals den Begriff Origami, um  eine Konstruktion W. Thurstons in \cite{thur} und deren Verallgemeinerung durch W. Veech in \cite{veech} zu beschreiben. Dabei sei bemerkt, dass es sich bei diesen mathematischen Objekten keinesfalls um die gleichnamige japanische Papierfaltkunst handelt. 

\newpage
{\bf Beispiel 1.1.2} 
\begin{list}{\alph{enumi})}{\usecounter{enumi} \setlength{\topsep}{0.2cm} \setlength{\parskip}{0.4cm}}
\item Das einfachste Beispiel ist ein Origami, welches nur aus einer Kopie des euklidischen Einheitsquadrates besteht. Denn daf"ur gibt es genau eine Verklebem"oglichkeit, die den obigen Regeln entspricht. Man erh"alt einen Torus $T$. Dieses Orgami hei"st triviales Origami $O_{0}$  und hat Geschlecht g=1.
\begin{center}
\begin{picture}(1,-2) 
\multiput(0,0)(0,-1){2}{\line (1,0){1}}
\multiput(0,0)(1,0){2}{\line (0,-1){1}}
\put(-0.1,-0.1){$\bullet$}
\put(0.9,-0.1){$\bullet$}
\put(-0.1,-1.1){$\bullet$}
\put(0.9,-1.1){$\bullet$} 
\end{picture} 
\end{center} 
\vspace{0.8cm}
\begin{center}
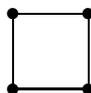
{\small Abb. 1: Das triviale Origami $O_{0}$. Gegen"uberliegende Kanten werden verklebt.}
\end{center}
\item Das Origami bestehend aus zwei Kopien des euklidischen Einheitsquadrates mit den in Abbildung 2 angegebenen Verklebungen hei"st $O_{1}$ und hat ebenfalls Geschlecht g=1.\begin{center}
\begin{picture}(2,-2)
\multiput(0,0)(1,0){2}{\line (1,0){1}}
\multiput(0,0)(1,0){3}{\line (0,-1){1}}
\multiput(0,-1)(1,0){2}{\line (1,0){1}}
\put(0.4,-0.6){$1$}
\put(1.4,-0.6){$2$}
\put(-0.1,-0.1){$\bullet$}
\put(0.9,-0.1){$\circ$}
\put(1.9,-0.1){$\bullet$}
\put(-0.1,-1.1){$\circ$}
\put(0.9,-1.1){$\bullet$}
\put(1.9,-1.1){$\circ$}
\put(0.4,0.1){$a$}
\put(0.4,-1.3){$b$}
\put(1.4,0.1){$b$}
\put(1.4,-1.3){$a$}
\end{picture}
\end{center}
\vspace{1cm}
\begin{center}
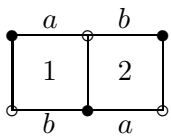
{\small Abb.2: Das Origami $O_{1}$. Kanten mit der selben Beschriftung und gegen"uberliegende Kanten werden verklebt.}
\end{center}
\item Das Origami bestehend aus f"unf Kopien des euklidischen Einheitsquadrates mit den in Abbildung 3 angegebenen Verklebungsvorschriften hei"st $D$. Man erh"alt in diesem Fall drei Identifikationsklassen von Ecken. Mit Hilfe der Eulerschen Charakteristik berechnet man f"ur die enstandene Fl"ache $X$ das Geschlecht g=2.
\begin{center}
\begin{picture}(3,3)
\multiput(0,0)(1,0){3}{\line (1,0){1}}
\multiput(0,0)(1,0){4}{\line (0,-1){1}}
\multiput(0,-1)(1,0){3}{\line (1,0){1}}
\multiput(0,0)(1,0){2}{\line (0,1){1}}
\multiput(0,1)(1,0){2}{\line (0,1){1}}
\multiput(0,1)(0,1){2}{\line (1,0){1}}
\put(0.4,-0.6){$1$}
\put(1.4,-0.6){$2$}
\put(2.4,-0.6){$3$}
\put(0.4,0.4){$4$}
\put(0.4,1.4){$5$}
\multiput(-0.1,-0.1)(1,0){2}{$\ast$}
\put(1.9,-0.1){$\bullet$}
\put(2.9,-0.1){$\ast$}
\multiput(-0.1,-1.1)(1,0){2}{$\bullet$}
\put(1.9,-1.1){$\ast$}
\put(2.9,-1.1){$\bullet$}
\multiput(-0.1,0.9)(1,0){2}{$\circ$}
\multiput(-0.1,1.9)(1,0){2}{$\bullet$}
\put(2.4,0.1){$a$}
\put(2.4,-1.3){$b$}
\put(1.4,0.1){$b$}
\put(1.4,-1.3){$a$}
\end{picture}
\end{center}
\vspace{1cm}
\begin{center}
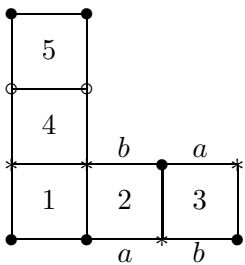
{\small Abb.3: Das Origami $D$. Kanten mit der selben Beschriftung und gegen"uberliegende Kanten werden verklebt.}
\end{center}
\end{list}	
\newpage
Weitere M"oglichkeiten, Origamis zu beschreiben, findet man in \cite{gabi1}. Dort werden unter anderem die folgenden Beschreibungen eingef"uhrt: 
\begin{list}{$\bullet$}{ \setlength{\topsep}{0cm} \setlength{\parskip}{0.4cm}}
\item Als "Uberlagerungen des Torus, die h"ochstens "uber einem Punkt verzweigen.
\item Als Paare von Permutationen in {\it$S_n$}, die gewisse Bedingungen erf"ullen.
\item Als Untergruppen der freien Gruppe $\mathbb{F}_2$ mit endlichem Index. 
\end{list}

F"ur die vorliegende Arbeit wird vor allem die in Punkt zwei angesprochene Definition verwendet. Die weiteren werden sp"ater in Kapitel 4 noch einmal aufgegriffen. 

{\bf Definition 1.1.3}\\
Ein {\it Origami O} bestehend aus $n$ Kopien des euklidischen Einheitsquadrates wird definiert durch zwei Permuationen $\sigma_a$ und $\sigma_b$ in {\it$S_n$}, wobei $\sigma_a$ und $\sigma_b$ induzieren, wie die vertikalen bzw. horizontalen Kanten verklebt werden. 

\textsc{Bermerkung:}\\
W"ahlt man zwei beliebige Permutationen $\sigma_a$ und $\sigma_b \in S_n$, so zeigt das folgenden Lemma, dass das Erzeugnis $\left\langle \sigma_a, \sigma_b\right\rangle$ der Permutationen transitiv auf $n$ operiert. Diese Voraussetzung ist n"otig, um eine zusammenh"angend Fl"ache zu erhalten.

{\bf Lemma 1.1.4}\\ \begin{it}
F"ur $n\rightarrow\infty$ und $\sigma_a, \sigma_b \in S_n$ operiert $\left\langle \sigma_a,\sigma_b \right\rangle$ transitiv.\end{it} \\
\bew Es gilt, dass zwei beliebig gew"ahlte Elemente der alternierenden Gruppe {\it$A_n$} diese mit Wahrscheinlichkeit $\rightarrow 1$ erzeugen, falls $n\rightarrow\infty$. Dies folgt aus Dixons Vermutung von 1969 und deren Beweis durch Dixon selbst \cite{dixon}, Kantor und Lubotzky \cite{kantlub} und Liebeck und Shalev \cite{liebshal}. Weiterhin gilt nach \cite{MR0183775}, dass {\it$A_n$} $(n-2)$-fach transitiv ist, das hei"st es existiert zu je zwei Folgen $a_1,...,a_{n-2}$ und $b_1,...,b_{n-2}$ aus $\left\{1,2,...,n\right\}$ von je $n-2$ verschiedenen Elementen eine Permutation $\sigma$  $\in$ {\it$A_n$}, so dass $\sigma(a_1)=b_1,..., \sigma(a_{n-2})=b_{n-2}$. Insbesondere ist {\it$A_n$}, f"ur $n\geq 3$, einfach transitiv, woraus die Aussage des Lemmas folgt. \qed

{\bf Definition 1.1.5}\\
Zwei Origamis $O$ und $O'$ hei"sen {\it "aquivalent}, falls die beiden Paare $(\sigma_a,\sigma_b)$ und $(\sigma_a',\sigma_b')$ in {\it$S_n$} "aquivalent sind. Das hei"st, falls ein $s\in$ {\it$S_n$} existiert, so dass $\sigma_a = s\sigma_a's^{-1}$ und $\sigma_b = s\sigma_b's^{-1}$. 

{\bf Beispiel 1.1.6}\\
In Beispiel 1.1.2 b) erh"alt man f"ur das Origami $O_1$ die Permutationen $\sigma_a = (1 2)$ und $\sigma_b = (1 2)$ in {\it$S_2$}. \\
F"ur das Origami $D$ in c) erh"alt man die Permuationen $\sigma_a = (1 2 3)$ und $\sigma_b = (1 4 5)(2 3)$ in {\it$S_5$}. Die anschauliche Bedeutung von $\sigma_a$ ist, dass \\
\hspace*{0.5cm}1. die rechte Kante des ersten Quadrates mit der linken des zweiten Quadrates, \\
\hspace*{0.5cm}2. die rechte Kante des zweiten Quadrates mit der linken des dritten Quadrates und \\
\hspace*{0.5cm}3. die rechte Kante des dritten Quadrates mit der linken des ersten Quardrates \\ identifiziert werden soll.\\
Der Zykel $(1 4 5)$ in $\sigma_b$ bedeutet, dass\\
\hspace*{0.5cm}1. die obere Kante des ersten Quadrates mit der unteren des vierten Quadrates, \\
\hspace*{0.5cm}2. die obere Kante des vierten Quadrates mit der unteren des f"unften Quadrates und \\
\hspace*{0.5cm}3. die obere Kante des f"unften Quadrates mit der unteren des ersten Quardrates \\
identifiziert werden soll.\\
Und der Zykel $(2 3)$ in $\sigma_b$ bedeutet, dass\\
\hspace*{0.5cm}1. die obere Kante des zweiten Quadrates mit der unteren des dritten Quadrates und \\
\hspace*{0.5cm}2. die obere Kante des dritten Quadrates mit der unteren des zweiten Quadrates\\ identifiziert werden soll.
\section{Kranzprodukt}
{\bf Definition 1.2.1}\\
Sei \hfill $G$ \hfill eine \hfill Gruppe \hfill und \hfill $H$ \hfill eine \hfill Permutationsgruppe \hfill auf \hfill der \hfill Menge \hfill der \hfill Symbole \\ $\Omega=\left\{1,...,n\right\}$, so hei"st die Menge \[\left\{(h;f) | h \in H, f:\Omega \rightarrow G\right\}\] zusammen mit der Verkn"upfungsvorschrift \[(h;f)(h';f'):=(hh';ff'_h)\] {\it Kranzprodukt $H\wr G$ von H mit G}. \\ F"ur $f:\Omega \rightarrow G$ und $h \in$ H ist die Abbildung $f_h:\Omega \rightarrow G$ definiert durch \[f_h(i):=f(h^{-1}(i)), \forall i \in \Omega\] und f"ur zwei Abbildungen $f,f':\Omega \rightarrow G$ ist deren Produkt $ff':\Omega \rightarrow G$ definiert durch \[ff'(i):=f(i)f'(i), \forall i \in \Omega.\] 
Sei $e:\Omega \rightarrow G$ die Abbildung mit den Werten $e(i)=1_G, \forall i \in \Omega$ und sei $f^{-1}:\Omega \rightarrow G$ definiert durch $f^{-1}:=f(i)^{-1}, \forall i \in \Omega$, so erh"alt man f"ur das Einselement in $H\wr G$ und f"ur das Inverse von $(h;f)$:
\[1_{H\wr G}=(1_H;e),\]
\[(h;f)^{-1}=(h^{-1};f^{-1}_{h^{-1}}),\]
wobei $f^{-1}_{h^{-1}}:=(f_{h^{-1}})^{-1}$ ist .\\
Es ist leicht zu zeigen, dass $H\wr G$ eine Gruppe der Ordnung $|H\wr G|=|G|^n|H|$ ist.
 
\textsc{Notation:}\\
F"ur $(h;f) \in$ $H\wr G$ kann ebenso die entsprechende Tabellenform $(h;f_1,\ldots ,f_n)$ verwendet werden, wobei $f_i:=f(i)$ ist.

{\bf Lemma 1.2.2}\\ \begin{it}
Sei $G\ni g \rightarrow \binom{i}{g(i)}$ eine Permutationsdarstellung von $G$ auf der Menge der Symbole $\Gamma=\left\{1,...,m\right\}$, so erh"alt man durch 
\begin{eqnarray} \label{kranzperm}
(h;f)(i,j):=(h(i);f_{h(i)}(j)), \forall (i,j) \in \Omega\times\Gamma
\end{eqnarray} eine Permutationsdarstellung von $H\wr G$ auf $\Omega\times\Gamma$.\end{it} \\
\bew Es ist einfach zu "uberpr"ufen, dass \[(h;f)((h';f)(i;j))=(hh';ff'_h)(i,j)\ \forall f,f',h,h',i,j\] gilt und aus $(h;f)(i,j)=(h;f)(i',j')$ folgt mit Gleichung (\ref{kranzperm}) \[(h(i);f_{h(i)}(j))=(h(i');f_{h(i')}(j')).\]
Da $h$ eine Permutation ist, folgt $i=i'$ und man erh"alt $f_{h(i)}(j)=f_{h(i)}(j')$, woraus wiederum, da $f_{h(j)}$ eine Permutation ist, folgt, dass $j=j'$.\qed

Im Folgenden sollen Kranzprodukte der Form $S_n\wr G$ betrachtet werden, wobei $G$ endlich sein soll.

{\bf Definition 1.2.3}\\
Sei $(h;f)$ $\in S_n\wr G$ und sei \[h=\prod^{c(h)}_{\nu=1} h_\nu=\prod^{c(h)}_{\nu=1}(j_\nu h(j_\nu)\cdots h^{l_{\nu-1}}(j_{\nu}))\] die Zykelnotation von $h$, wobei $c(h)$ die Anzahl der disjunkten zyklischen Faktoren einschlie"slich der 1-Zykel ist, $l_\nu$ deren L"ange und $j_\nu$ das kleinste Symbol im zyklischen Faktor $\nu$. So definiert
\begin{eqnarray} \label{zykelprodukt}
g_\nu(h;f) := f(j_\nu)f(h^{-1}(j_\nu))\cdots f(h^{-l_\nu+1}(j_\nu)) = ff_h\cdots f_{h^{l_\nu-1}}(j_\nu) 
\end{eqnarray}
das {\it $\nu$-te Zykelprodukt von $(h;f)$}.

{\bf Definition 1.2.4}\\
Sei die Permutation $h$ $\in S_n$ vom Zykeltyp $Th=(a_1,\ldots ,a_n)=(1^{a_1},\ldots ,n^{a_n})$, das hei"st $h$ besitzt $a_k$ zyklische Faktoren der L"ange $k$, und sei $f$ eine Abbildung von $\Omega$ nach $G$, so dass $(h;f)$ $\in S_n\wr G$. 
So gibt es $a_k$ zu den $a_k$ Zykeln der L"ange k geh"orige Zykelprodukte. Seien $C^1,\ldots ,C^s$ die Konjugationsklassen von $G$ und geh"oren $a_{ik}$ der $a_k$ Zykelprodukte zu $C^i$, so ist 
\begin{eqnarray*}
T(h;f):=(a_{ik})_{1\leq i \leq s \atop 1 \leq k \leq n}
\end{eqnarray*}
der {\it Zykeltyp von $(h;f)$}.

Es lassen sich nun die Konjugationsklassen von $S_n\wr G$ beschreiben. 

{\bf Lemma 1.2.5} \begin{it} \\
Zwei \hfill Elemente \hfill $(h;f)$ \hfill und \hfill $(h';f')$ \hfill von \hfill $S_n\wr G$ \hfill sind \hfill konjugiert \hfill genau \hfill dann, \hfill wenn \\ $T(h;f)=T(h';f')$ gilt.\end{it} \\
\bew siehe \cite[S.141, Thm 4.2.8]{jk}. \qed 
\chapter{Darstellungstheorie}
\section{Darstellungstheorie endlicher Gruppen}
Dieses Kapitel soll die f"ur diese Arbeit relevanten Ergebnisse der Charaktertheorie endlicher Gruppe vorstellen. Hierzu werden zun"achst einige grundlegende Definitionen und Resultate der Darstellungstheorie angegeben.

{\bf Definition 2.1.1}\\
Sei $V$ ein Vektorraum "uber einem K"orper $K$, $GL(V)$ die Gruppe aller Automorphismen von $V$ und $G$ sei endliche Gruppe, so ist eine {\it Darstellung von G "uber K} ein Paar $(D,V)$ mit $D_V: G \rightarrow Gl(V)$ Gruppenhomomorphismus. Der Grad der Darstellung ist die Vektorraumdimension $n$. 

\textsc{Notation:}\\
Ist klar, welche Abbildung gemeint ist, wird $V$ meist selbst als Darstellung (oder $K$-Darstellung) von $G$ bezeichnet. F"ur $D(g)(v)$ schreibt man kurz $g.v$ f"ur $g\in G, v\in  V$. 

\textsc{Bemerkung:}\\
Diese Schreibweise macht deutlich, dass man jede Darstellung von $G$ auch als links $G$-Modul auffassen kann, da jede Darstellung $D$ von $G$ eine Operation von $G$ auf dem $K$-Vekorraum $V$ definiert.
\newpage
\textsc{Bemerkung:}\\
Sei $(D,V)$ Darstellung von $G$ und $H$ Untergruppe von $G$, so induziert die auf Elemente aus $H$ eingeschr"ankte Darstellung $(D,V)$  eine Darstellung $(D_{|H},V)$ von $H$. 

{\bf Beispiel 2.1.2} \begin{list}{\alph{enumi})}{\usecounter{enumi} \setlength{\topsep}{0cm} \setlength{\parskip}{0.4cm}}
\item Die triviale Darstellung bildet alle Elemente $g\in G$ auf die Identit"at ab.
\item Sei $m$ die Ordnung von $G$ und sei $V$ ein Vektorraum der Dimension $m$ mit Basis $(e_t)_{t\in G}$. F"ur $g\in G$ sei $D_g$ die lineare Abbildung von $V$ nach $V$, welche $e_t$ nach $e_{gt}$ abbildet. Die Vorschrift $g \mapsto D_g$ definert die sogenannte regul"are Darstellung.
\end{list}
{\bf Definition 2.1.3}\\
Seien $(D,V)$ und $(D',V')$ zwei Darstellungen von $G$, so hei"sen diese {\it isomorph oder "aquivalent}, falls ein Isomorphismus $f: V \rightarrow V'$ existiert mit \[f(D_V(g)(v))=D'_{V'}(g)(f(v)) \ \forall g \in G, v \in V.\]

{\bf Definition 2.1.4}\\
Sei $G$ endliche Gruppe, $H$ Normalteiler von $G$ und $(D_{|H},V)$ Darstellung von $H$. $\mathbb{D}(h)$ sei die zu $D_{|H}(h)$ korrespondierende Matrix und $g \in G$ fix, so hei"st die durch 
\[h\mapsto \mathbb{D}(g^{-1}hg)\] definierte Darstellung von $H$, die zu $D_{|H}$ {\it konjugierte Darstellung} bez"uglich $G$.

{\bf Definition 2.1.5}\\
a) Sei $V$ Darstellung von $G$ und $U\subseteq V$ Untervektorraum, so hei"st $U$ {\it Unterdarstellung von V}, falls $g.u \in U$  $\forall g \in G, u \in U$.\\
b) Eine Darstellung $V\neq\left\{0\right\}$ hei"st {\it irreduzibel}, falls $\left\{0\right\}$ und $V$ die einzigen Unterdarstellungen sind.

{\bf Definition 2.1.6}\\
Eine Darstellung $(D,V)$ von $G$ hei"st {\it induziert} durch eine Darstellung $(R,U)$ von $H$, falls $V={\bigoplus \atop \scriptstyle \sigma \in G/H} \sigma U$, wobei $\sigma U:=D_gU \subset V$ f"ur ein $g \in \sigma$.

{\bf Satz 2.1.7}\\ \begin{it}
Sei $(R,U)$ Darstellung von $H$, so exisiert eine Darstellung $(D,V)$ von $G$, welche durch $(R,U)$ induziert wird. Diese Darstellung ist eindeutig bis auf Isomorhismen.\end{it}\\
\bew  siehe \cite[S.29f]{se}. \qed 

{\bf Definition 2.1.8}\\
Seien $V,V'$ Darstellungen von $G$ so ist $V\oplus V'$ mittels $g.(v,v') = (g.v,g.v')$ wieder eine Darstellung von $G$. $V\oplus V'$ hei"st {\it direkte Summe von V und V'}. 

{\bf Definition 2.1.9}\\
Sei $V$ endlich dimensionale Darstellung von G "uber dem K"orper der komplexen Zahlen $\mathbb{C}$. Der {\it Charakter $\chi$ von V} ist  
\begin{eqnarray*}
\chi=\chi_V : & G & \rightarrow K\\
& g & \longmapsto Spur(\mathbb{D}(g)).
\end{eqnarray*} 

{\bf Lemma 2.1.10} \begin{it} \\
Sei $\chi$ Charakter einer Darstellung vom Grad $n$ und $1\in G$ Einheit, so gilt:\\
\hspace*{0.5cm}$(i)$ $\chi(1)=n$, \\
\hspace*{0.5cm}$(ii)$ $\chi(g^{-1})=\overline{\chi(g)}\ \forall g \in G$,\\
\hspace*{0.5cm}$(iii)$ $\chi(hgh^{-1})=\chi(g)\ \forall g \in G$,\\
\hspace*{0.5cm}$(iv)$ $\chi_{V\oplus W}=\chi_V+\chi_W$.\\
\end{it}
\bew $(i)$ $\chi(1)= Spur(id_V)=dim V=n$\\
$(ii)$ F"ur alle $g \in G$ existiert, da $G$ endlich ist, ein $N\ge 0$, so dass $g^N=1$.\\ $\Rightarrow$ Alle Eigenwerte $\lambda_1,\ldots ,\lambda_n$ von $\mathbb{D}(g)$ haben Betrag 1 und somit \[ \overline{\chi(g)}=\overline{Spur(\mathbb{D}(g))}=\sum {\overline{\lambda_i}}=\sum {\lambda_i^{-1}}= Spur(\mathbb{D}(g)^{-1})= Spur (\mathbb{D}(g^{-1})) = \chi(g^{-1})\]
$(iii)$ $\chi(hgh^{-1}) = Spur(\mathbb{D}(hg)\mathbb{D}(h^{-1})) = Spur(\mathbb{D}(h^{-1})\mathbb{D}(hg)) = \chi(g)$.\\
$(iv)$ $\mathbb{D}_{V\oplus W}(g)=\begin{pmatrix} \mathbb{D}_V(g) & 0 \\ 0 & \mathbb{D}_W(g) \end{pmatrix}$\\
 so dass gilt $Spur(\mathbb{D}_{V\oplus W}(g))=Spur(\mathbb{D}_V(g))+Spur(\mathbb{D}_W(g))$. \qed 
 
{\bf Satz (Orthogonalit"atsrelationen) 2.1.11} \begin{it} \\
Seien $V,W$ zwei irreduzible Darstellungen, $\chi_V, \chi_W$ deren Charaktere und sei \[(\chi_V,\chi_W):= \frac{1}{|G|} \sum_{g\in G} \chi_V(g)\chi_W(g) \in \mathbb{C},\] so gilt:\\
\hspace*{0.5cm} a) $(\chi_V,\chi_V)=1$\\
\hspace*{0.5cm} b) $V\not\cong W$ $\Rightarrow$ $(\chi_V, \chi_W)=0$\\
\end{it} 
\bew siehe \cite[S. 15, Thm 3]{se}. \qed 

Das folgende Theorem besagt, wie man den Charakterwert einer induzierten Darstellung berechnet. Sei dazu wie oben $(D,V)$ durch $(R,U)$ induziert und $\chi^D$ und $\chi^R$ seien die zugeh"origen Charaktere von $G$ bzw. $H$.

{\bf Theorem 2.1.12} \begin{it}\\
Sei $h$ die Ordnung von $H$ und sei $S$ ein Repr"asentantensystem von $G/H$. So l"asst sich der Charakterwert f"ur $g\in G$ durch
\bq \label{ind}
\chi^D(g)=\sum_{s\in S \atop sgs^{-1} \in H} \chi^R(sgs^{-1})=\frac{1}{h}\sum_{g'\in G \atop g'gg'^{-1} \in H} \chi^R(g'gg'^{-1})
\eq 
berechnen.\end{it}\\
\bew siehe \cite[S.30, Thm 12]{se}. \qed 

{\bf Satz 2.1.13} \begin{it} \\
Sei V eine endlich dimensionale Darstellung  von $G$, $W$ irreduzibel und $V\cong W_1\oplus \cdots \oplus W_n$ mit $W_i$ irreduzibel. Dann ist die Anzahl $m_w=\#\left\{i=1,\ldots ,n | W_i\cong W\right\}=(\chi_V,\chi_W)$. Insbesondere ist $m_w$ unabh"angig von der Zerlegung von $V$.\end{it} \\
\bew $\chi_V=\sum_{i=1}^n \chi_{W_i}$, $(\chi_V, \chi_W)=\sum_{i=1}^n (\chi_{W_i}, \chi_W)$ mit Satz 2.1.8 folgt nun die Aussage des Satzes, da $(\chi_{W_i}, \chi_W)$ entweder 0 oder 1 ist, abh"angig davon ob $W_i$ isomorph zu $W$ ist oder nicht. \qed 

{\bf Lemma 2.1.14} \begin{it} \\
Zwei endlich dimensionale Darstellungen sind genau dann isomorph, wenn ihre Charaktere gleich sind.\end{it}\\
\bew Satz 2.1.13 zeigt, dass eine irreduzible Darstellung in jeder der beiden Darstellungen mit der selben Multiplizit"at vorkommt, falls die Charaktere "ubereinstimmen. Daher sind in diesem Fall die beiden Darstellungen isomorph. Ebenso gilt, dass zwei isomorphe Darstellungen diesselbe Spur besitzen, so dass ihre Charaktere "ubereinstimmen.\qed 

{\bf Definition 2.1.15}\\
Eine Abbildung $f: G \rightarrow \mathbb{C}$ hei"st {\it Klassenfunktion}, falls sie konstant auf den Konjugationsklassen ist.

\textsc{Bemerkung:}\\
$H$:=$\left\{f: G \rightarrow \mathbb{C} | \mbox{f ist Klassenfunktion} \right\}$ ist ein endlich dimensionaler $\mathbb{C}$-Vektorraum. 

{\bf Satz 2.1.16} \begin{it} \\
Seien $W_1,\ldots ,W_n$ irreduzible Darstellungen von $G$ und $\chi_1,\ldots ,\chi_n$ ihre Charaktere, dann ist $\chi_1,\cdots ,\chi_n$ eine Orthonormalbasis von $H$.\end{it} \\
\bew siehe \cite[S.19, Thm 6]{se}. \qed 
\newpage
{\bf Lemma 2.1.17} \begin{it} \\
Die Anzahl der irreduziblen Darstellungen von $G$ (bis auf Isomorphie) ist gleich der Anzahl der Konjugationsklassen.\end{it} \\
\bew Seien $C_1,\ldots C_k$ die Konjugationsklassen von $G$, so ist eine Klassenfunktion $f$ durch die Werte $\lambda_i$ auf den Klassen $C_i$ eindeutig bestimmt. Also ist die Dimension von $H$ gleich k und nach Satz 2.1.16 ist die Dimension gerade die Anzahl der irreduziblen Darstellungen von $G$ (bis auf Isomorphie). Mit Lemma 2.1.14 folgt dann die Aussage des Lemmas. \qed 

{\bf Lemma 2.1.18} \begin{it} \\
Sei $G$ eine endliche Gruppe und seien $C_1,C_2$ zwei Konjugationsklassen von $G$ mit Repr"asentanten $g_i$. So ist die Anzahl der L"osungen der Gleichung $g_1g_2=z$ gleich 
\begin{eqnarray} \label{frob}
\frac{|C_1||C_2|}{|G|} \sum_{\chi \in Irr(G)} \frac{\chi(g_1)\chi(g_2)\chi(z^{-1})}{\chi(1)},
\end{eqnarray}
 wobei $Irr(G)$ die Menge aller irreduziblen Darstellungen von $G$ ist.\end{it} \\
\bew siehe \cite[S.216, Prop. 9.33]{cu}.\qed
\section{Darstellungstheorie der symmetrischen Gruppe}
Diese Theorie wurde gr"o"stenteils durch die Entdeckung der Young Tableaus durch Alfred Young in dessen Arbeiten "uber die Invariantentheorie \cite{y1} vorangetrieben, da diese Tableaus die irreduziblen Darstellungen von $S_n$ beschreiben.\\
Dieses Kapitel beinhaltet wiederum nur die f"ur diese Arbeit relevanten Bereiche der Charaktertheorie der symmetrischen Gruppe. Ausf"uhrlichere Beschreibungen findet man in \cite{rob} und \cite{jk}.

\textsc{Bemerkung:}\\
Die Konjugationsklassen der symmetrischen Gruppe $S_n$ sind definiert durch eine Partition $(\lambda)=(\lambda_1,\cdots \lambda_h)$, wobei $\lambda_1+\ldots +\lambda_h=n$. Die Norm $\|\lambda\|$ von $\lambda$ ist das Gr"o"ste $j$, so dass $\lambda_j\neq 0$. 
\newpage
{\bf Definition 2.2.1} \\
a) Eine Partition ($\lambda$) kann durch ein korrespondierendes {\it Young Diagramm} $[\lambda]$, welches \hspace*{0.4cm} aus $n$ Punkten besteht, dargestellt werden. Die i-te Zeile von $[\lambda]$ besteht dabei aus $\lambda_i$ \hspace*{0.4cm} Punkten und jede Zeile beginnt in der selben Spalte: 
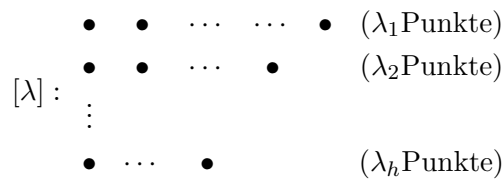
\[
\left[\lambda\right]:
\begin{array}{cccccc}
  \bullet & \bullet & \cdots & \cdots & \bullet & (\lambda_1 \mbox{Punkte}) \\
  \bullet & \bullet & \cdots & \bullet & &  (\lambda_2 \mbox{Punkte}) \\
  \vdots & & & & &\\
  \bullet & \cdots & \bullet & & &  (\lambda_h \mbox{Punkte}) \end{array} \] 
\hspace*{0.4cm} $\left[\lambda\right] \geq \left[\mu\right]$, falls $\lambda_1=\mu_1, \lambda_2=\mu_2,\ldots ,\lambda_r=\mu_r, \lambda_{r+1}\geq \mu_{r+1}$.\\  
b) Seien $[\lambda]$ und $[\mu]$ zwei Young Diagramme und sei $[\mu]$ vollst"andig in $[\lambda]$ enthalten, so \hspace*{0.4cm} hei"st das Residuum $[\lambda]- [\mu]$ {\it Skew Diagramm}. \\
c) Das zu $[\lambda]$ {\it konjugierte} Diagramm $[\lambda']$ erh"alt man durch Vertauschung der Zeilen und \hspace*{0.4cm} Spalten.\\
d) Ersetzt man im Young Diagramm die Punkte durch die Symbole 1,2,$\ldots ,n$, so gibt \hspace*{0.4cm} es $n!$ M"oglichkeiten, diese Symbole anzuordnen. Diese Diagramme hei"sen dann {\it Young \hspace*{0.4cm} Tableaus} und werden mit $t_1^\lambda,t_2^\lambda,\ldots ,t_{n!}^\lambda$ bezeichnet.\\
e) W"ahlt man aus diesen $n!$ Young-Tableaus diejenigen aus, deren Eintr"age aufsteigend \hspace*{0.4cm} in Zeile und Spalte sind, so nennt man diese {\it Standard Tableaus}. Deren Anzahl wird \hspace*{0.4cm} mit $f^\lambda$ bezeichnet.\\
\hspace*{0.4cm} Auf diesen Standard Tableaus l"asst sich eine {\it lexikographische Ordnung} definieren:\\
\hspace*{0.4cm} $t_i^\lambda \geq t_j^\lambda$, falls die Symbole der ersten r Zeilen und die ersten s Symbole der (r+1)-ten \hspace*{0.4cm} Zeile beider Tableaus "ubereinstimmen, aber das (s+1)-te Symbol von $t_i^\lambda$ gr"o"ser ist als \hspace*{0.4cm} das (s+1)-te Symbol von $t_j^\lambda$.
 
Nach Lemma 2.1.17 werden die irreduziblen Darstellungen der symmetrischen Gruppe $S_n$ nun mit den korrespondierenden Young Tableaus $[\lambda]$ und die dazugeh"origen Charaktere mit $\chi^\lambda$ bezeichnet.

{\bf Definition 2.2.2} \\
a) Der $Haken \ (i,j)$ von $[\lambda]$ besteht aus dem Punkt $(i,j)$ zusammen mit den ($\lambda_i-i$) \hspace*{0.4cm} Punkten rechts davon, dem {\it Arm}, und den $(\lambda_j'-j)$ Punkten darunter, dem {\it Bein}.\\
b) Die {\it L"ange} des Hakens $(i,j)$ ist $h_{ij}=(\lambda_i-i)+(\lambda_j'-j)+1$.\\
c) Ersetzt man das Symbol an der Stelle $(i,j)$ durch $h_{ij}$, so erh"alt man das {\it Haken Tableau} \hspace*{0.4cm} $H[\lambda]= (h_{ij})$. Mit $H^\lambda$ bezeichnet man das Produkt aller $h_{ij}$.\\
d) Ein {\it Skew-Haken} ist ein zusammenh"angender Teil des Randes eines Young Diagramms, \hspace*{0.4cm} so dass bei dessen Entfernung wieder ein Young Diagramm entsteht. Die L"ange eines \hspace*{0.4cm} Skew-Hakens entspricht der Anzahl der enthaltenen Punkte. 

{\bf Beispiel 2.2.3} \\
Das zur Partition $\lambda=[4^2,3]$ geh"orende Young-Diagramm besitzt 
\begin{list}{$\bullet$}{\setlength{\topsep}{0cm} \setlength{\parskip}{0.4cm}}
\item ein Skew-Haken der L"ange 6:
\begin{center}
\begin{picture}(1.5,-1)
\multiput(0,0)(0.5,0){4}{$\bullet$}
\multiput(0,-0.5)(0.5,0){4}{$\bullet$} 
\multiput(0,-1)(0.5,0){3}{$\bullet$}
\put(1.1,-0.4){\line(1,0){0.5}}
\put(1.6,-0.4){\line(0,1){0.5}}
\put(1.1,-0.9){\line(0,1){0.5}}
\multiput(0.1,-0.9)(0.5,0){2}{\line(1,0){0.5}}
\end{picture}
\end{center}
\vspace{0.7 cm}
\item zwei Skew-Haken der L"ange 5:
\begin{center}
\begin{picture}(6,-1)
\multiput(0,0)(0.5,0){4}{$\bullet$}
\multiput(0,-0.5)(0.5,0){4}{$\bullet$} 
\multiput(0,-1)(0.5,0){3}{$\bullet$}
\put(1.1,-0.4){\line(1,0){0.5}}
\put(1.1,-0.4){\line(0,-1){0.5}}
\multiput(0.1,-0.9)(0.5,0){2}{\line(1,0){0.5}}
\multiput(4.5,0)(0.5,0){4}{$\bullet$}
\multiput(4.5,-0.5)(0.5,0){4}{$\bullet$} 
\multiput(4.5,-1)(0.5,0){3}{$\bullet$}
\put(5.6,-0.4){\line(1,0){0.5}}
\put(5.6,-0.4){\line(0,-1){0.5}}
\put(6.1,-0.4){\line(0,1){0.5}}
\put(5.6,-0.9){\line(-1,0){0.5}}
\end{picture}
\end{center}
\vspace{0.7 cm}
\item zwei Skew-Haken der L"ange 4:
\begin{center}
\begin{picture}(6,-1)
\multiput(0,0)(0.5,0){4}{$\bullet$}
\multiput(0,-0.5)(0.5,0){4}{$\bullet$} 
\multiput(0,-1)(0.5,0){3}{$\bullet$}
\put(1.1,-0.4){\line(1,0){0.5}}
\put(1.1,-0.4){\line(0,-1){0.5}}
\put(0.6,-0.9){\line(1,0){0.5}}
\multiput(4.5,0)(0.5,0){4}{$\bullet$}
\multiput(4.5,-0.5)(0.5,0){4}{$\bullet$} 
\multiput(4.5,-1)(0.5,0){3}{$\bullet$}
\put(5.6,-0.4){\line(1,0){0.5}}
\put(5.6,-0.4){\line(0,-1){0.5}}
\put(6.1,-0.4){\line(0,1){0.5}}
\end{picture}
\end{center}
\vspace{0.7 cm}
\item zwei Haken der L"ange 3
\item zwei Skew-Haken der L"ange 2
\item und zwei Skew-Haken der L"ange 1.
\end{list}

{\bf Lemma 2.2.4} \begin{it} \\
Zwischen den Haken und den Skew-Haken von $[\lambda]$ gibt es eine nat"urliche 1-1 Korrespondenz. \\ \end{it}
\bew 
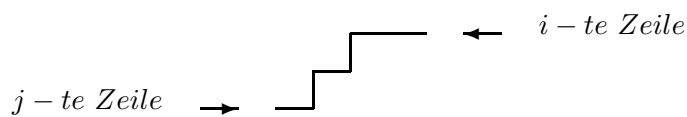
\begin{picture} (2,0)
\thicklines
\put(5,0){\line(1,0){1}}
\put(7,0){\vector(-1,0){0.5}}
\put(3,-1){\vector(1,0){0.5}}
\put(7.5,0){$ i-te$ $Zeile$}
\put(0.5,-1){$ j-te$ $Zeile$}
\put(5,0){\line(0,-1){0.5}}
\put(5,-0.5){\line(-1,0){0.5}}
\put(4.5,-0.5){\line(0,-1){0.5}}
\put(4.5,-1){\line(-1,0){0.5}}
\end{picture} 
\vspace{1cm}

Dieser Skew-Haken korrespondiert zu dem Haken $(i,j)$.\qed 

{\bf Korollar 2.2.5} \begin{it} \\
F"ur die Anzahl $f^\lambda$ der Standard Tableaus gilt:
\begin{eqnarray} \label{flambda}
f^\lambda= n!\frac{\prod\limits_{i < k}(h_{i1}-h_{k1})}{\prod\limits_i h_{i1}!}.
\end{eqnarray} \end{it} 
\bew siehe \cite[S.44 2.36]{rob}.\qed 

Das folgende Theorem gibt an, wie die zu einer Transposition korrespondierende Matrix explizit konstruiert werden kann. Im Gegensatz dazu folgt im Anschluss eine wesentlich einfachere Rekursionsformel f"ur die Berechnung der Charakterwerte belibieger Elemente $\pi \in S_n$.

{\bf Youngs Fundamental Theorem 2.2.6}\begin{it}\\
Die Konstruktion der Matrix $\mathbb{D}=(d_{ij})$, welche zur Transposition $(r,r+1)$ der irreduziblen Darstellung $[\lambda]$ korrespondiert, erfordert zun"achst eine lexikographische Anordnung der $f^\lambda$ Standard Tableaus $\ldots t_i^\lambda,\ldots ,t_j^\lambda,\ldots $. Setze dann
\begin{enumerate}
\item[(i)]
1 an die Stelle $d_{ii}$, falls $r$ und $r+1$ in der selben Zeile von $t_i^\lambda$ vorkommen,
\item[(ii)]
-1 an die Stelle $d_{ii}$, falls $r$ und $r+1$ in der selben Spalte von $t_i^\lambda$ vorkommen oder
\item[(iii)]
eine quadratische Matrix \\
$(a)$
$\bordermatrix{
  & t_i^\lambda & t_j^\lambda \cr
t_i^\lambda & -p & 1-p^2 \cr
t_j^\lambda & 1 & p \cr
}$  oder     
$(b)$
$\bordermatrix{
 & t_i^\lambda & t_j^\lambda \cr
t_i^\lambda & -p & \sqrt{1-p^2} \cr
t_j^\lambda & \sqrt{1-p^2} & p \cr
}$

an die Schnittstelle der zu $t_i^\lambda$ und $t_j^\lambda$ korrespondierenden Zeilen und Spalten, falls $i < j$ und $t_j^\lambda$ durch $t_i^\lambda$ ensteht, indem man $r$ und $r+1$ vertauscht. \\ 
Tritt $r$ an der Stelle $(k,l)$ und $r+1$ an der Stelle $(m,n)$ von $t_i^\lambda$ auf, wobei $k < m, \\ l > n$, so ist 
\[\frac{1}{p}=(l-k)(n-m).\]
\item[(iv)] 
0, sonst.
\end{enumerate} \end{it}
\bew siehe \cite[(VI) 34:196-230 (1932)]{y1}.\qed 

{\bf Beispiel 2.2.7}\\
Sei $\left[\lambda\right]=\left[2,2\right]$ Darstellung von $S_4$, so lassen sich die zu den Permutationen $(12)$ und $(12)(34)$ korresponierenden Matrizen nach Theorem 2.2.6 folgenderma"sen berechnen:
\begin{list}{\arabic{enumi})}{\usecounter{enumi} \setlength{\topsep}{0cm} \setlength{\parskip}{0.4cm}}
\item $f^\lambda=4! \frac{1!}{3! 2!}=2$ nach Gleichung (\ref{flambda})
\item $t_1^\lambda=\begin{array}{cc} 1 & 3 \\ 2 & 4 \end{array}$, $t_2^\lambda=\begin{array}{cc} 1 & 2  \\ 3 & 4 \end{array}$
\item Man erh"alt $\left[2,2\right]((12))=\begin{pmatrix}  -1 & 0 \\  0 & 1 \end{pmatrix}$ und $\left[2,2\right]((12)(34))=\begin{pmatrix} 1 & 0 \\  0 & 1 \end{pmatrix}.$
\end{list}
{\bf Murnaghan-Nakayama Formel 2.2.8} \begin{it} \\
Sei $\pi=\sigma\gamma$ disjunktes Produkt von $\sigma \in S_{n-k}$ und einem $k$-Zykel $\gamma$. So gilt:
\begin{eqnarray} \label{Mu-Na}
\chi^\lambda(\pi)=\sum_\mu \chi^{\lambda - \mu}(\gamma) \cdot \chi^\mu(\sigma),\end{eqnarray}
wobei die Summe "uber alle Diagramme $\left[\mu\right]$ bestehend aus $n-k$ Punkten l"auft. \end{it} \\
\bew siehe \cite[S.79]{james682}. \qed 

Die rekursive Berechnung der einzelen Charakterwerte von $S_n$ mit Hilfe obiger Formel, die erstmals von Murnaghan \cite{mur1} bewiesen wurde, l"asst sich noch durch das folgende Lemma vereinfachen.

{\bf Lemma 2.2.9} \begin{it} \\
Sei $[\lambda]-[\mu]$ Darstellung von $S_n$ und $\pi$ Zykel der L"ange $n$, so l"asst sich der Charakter von $\pi$ in der Darstellung $[\lambda]-[\mu]$ wie folgt berechnen: \begin{eqnarray} \label{lemma Mu-Na}
\chi^{\lambda-\mu}(\pi)=\begin{cases}
(-1)^l,& \mbox{falls} \ [\lambda]-[\mu] \ \mbox{ein Skew-Haken mit Beinl"ange l ist}\\
0, & \mbox{sonst}\\
\end{cases}\end{eqnarray}
wobei die Beinl"ange des zum Skew-Haken korrespondierenden Hakens gemeint ist. 
\end{it}\\
\bew siehe \cite[S.77 4.15]{rob}. \qed 

Die Vereinfachung der Formel ($\ref{Mu-Na}$) entsteht also dadurch, dass nicht mehr alle Diagramme $[\mu]$ auf den $n-k$ Punkten betrachtet werden m"usssen, sondern nur noch diejenigen, deren Residuum $[\lambda]-[\mu]$ ein Skew Haken ist.

{\bf Beispiel 2.2.10} 
\begin{enumerate}
\item[a)]
Im Fall von $S_2$ gibt es nur zwei irreduzible Charaktere $[2]$ und $[1^2]$. F"ur diese gilt nach Satz 2.1.10 $(i)$ und 2.1.11, dass $\chi^2(1^2)=\chi^{1^2}(1^2)=1$ und $\chi^2(2)=1$, $\chi^{1^2}(2)=-1$ ist.
\item[b)]
F"ur die Berechnung von $\chi^{2,1}(2)(1)$ gilt nach Gleichung (\ref{Mu-Na}) \[ \chi^{2,1}(2)(1)= 1\cdot \chi^2(2)+ 1\cdot \chi^{1^2}(2)= 1 \cdot 1+ 1 \cdot (-1)=0.\]
\item[c)]
$\chi^{3,1}(2)(2)$ l"asst sich mit Formel ($\ref{Mu-Na}$) und ($\ref{lemma Mu-Na}$) folgenderma"sen rekursiv berechnen \[
\chi^{3,1}(2)(2)=1 \cdot \chi^{3,1-2}(2)+ (-1) \cdot \chi^{3,1-1^2}(2)= 1 \cdot 0+ (-1) \cdot 1= -1\] 
\end{enumerate} 
Das folgende asymptotische Resultat der Charaktertheorie in Satz 2.2.11 und das anschlie"sende Theorem "uber die asymptotische Verteilung der Zykelanzahl einer Permutation von $n$ Elementen sind wichtige Hilfsmittel zur Berechnung der Wahrscheinlichkeit, dass das Produkt zweier beliebiger Elemente $\sigma, \tau \in S_n$ genau $s$ Zykel besitzt. 

{\bf Satz 2.2.11} \begin{it} \\
Es gilt: \begin{eqnarray} \label{lulov}
\sum_\lambda \chi^\lambda(1)^{-c}=2 \sum_{\lambda : \lambda_1 > n-A} \chi^\lambda(1)^{-c} + \mathcal{O}(n^{-Ac}).
\end{eqnarray} \end{it}
\bew siehe \cite{lulov}.\qed 

{\bf Theorem 2.2.12} \begin{it} \\
Sei $M_{(n)}$ die Anzahl der Zykel einer Permutation von $n$ Elementen, $E_n$ der Erwartungswert und $\sigma_n$ die Standardabweichung der Zufallsvariablen $M_{(n)}$, so erh"alt man:
\begin{enumerate}
\item[(i)] f"ur $\sigma_n$ und $E_n$ die folgenden Werte
\[E_n=\log n + \gamma + o(1)\]
\[\sigma_n = \sqrt{\log n} - (\frac{\pi^2}{12}-\frac{\gamma}{2})\frac {1}{\sqrt {\log n}} + o(\frac {1}{\sqrt {\log n}}),\] wobei $\gamma$ die Eulerkonstante ist 
\item[(ii)]
f"ur die charakteristische Funktion $\Phi_n(t)=\int_{-\infty}^\infty e^{itx} dF_n(x)$, wobei $F_n(x)$ Verteilungsfunktion der standardisierten Zufallsvariable $\frac{M_{(n)}-E_n}{\sigma_n \sqrt2}$:
\[ \lim_{n \to \infty} \Phi_n(t)= e^{-\frac{t^2}{4}}.\]
Also ist die standardisierte Zufallsvariable standardnormalverteilt.
\end{enumerate}
\end{it}
\bew siehe \cite[S.35 ff.]{gon}.\qed
\section{Darstellungstheorie des Kranzproduktes $H\wr G$}
Dieses Kapitel beschreibt in Anlehnung an \cite[Kapitel 4]{jk} und \cite[Kapitel 1]{kerber2} die Konstruktion irreduzibler Darstellungen von $H\wr G$ nach der Darstellungstheorie von A.H. Clifford in \cite{cliff}. Dabei sei $G$ eine endliche Gruppe und $H$ eine Permuatationsgruppe auf der Menge der Symbole $\Omega=\left\{1,...,n\right\}$ und $K$ im Folgenden stets ein algebraisch abgeschlossener K"orper. F"ur $K=\mathbb{C}$ hat W. Specht in \cite{specht2} erstmals die irreduziblen Darstellungen von $H\wr G$ konstruiert, ebenso wurde der Spezialfall $S_n\wr G$ von W. Specht in \cite{specht1} behandelt. 

{\bf Definition 2.3.1} \\
Die normale Untergruppe\[G^*:=\left\{(1_H;f)|f \in G^n\right\}\] hei"st {\it Basisgruppe} von $H\wr G$. Sie ist das direkte Produkt von $n$ Kopien von $G_i$, wobei \[G_i:=\left\{(1_H;f)|f_j=1_G \  \forall j\neq i\right\} \cong G.\] 
Die Untergruppe $H':=\left\{(h,e)| h\in H\right\}$, welche isomorph zu $H$ selbst ist, ist das {\it Komplement der Basisgruppe $G^*$}.

Sei $D$ eine irreduzible Darstellung von $G$ mit zugeh"origem Vektorraum $V$, so ist es Ziel dieses Kapitels zu zeigen, dass diese Darstellung zu einer Darstellung von $H\wr G$ erweitert werden kann. Sei dazu $\left\{b_1,\ldots ,b_m\right\}$ $K$-Basis von $V$, das hei"st $V= \left\langle \left\langle b_1,\ldots ,b_m \right\rangle\right\rangle_K$ und sei $\stackrel{n}{\bigotimes} V:=V\otimes_K \ldots \otimes_K V$ ($n$ Faktoren), so erh"alt man eine Basis von $ \stackrel{n}{\bigotimes} V$ durch:
\begin{eqnarray} \label{VR}
\bigotimes^n V=\left\langle \left\langle b_\varphi | \varphi \in \left[m\right]^{\left[n\right]}\right\rangle\right\rangle_K,
\end{eqnarray}
wobei $\left[m\right]=\left\{1,\ldots ,m\right\}$ bzw. $\left[n\right]=\left\{1,\ldots ,n\right\}$.\\
Fasst man $D$ nun als links $G$-Modul auf, so wird $\stackrel{n}{\bigotimes} V$ ein links $H\wr G$-Modul mittels der folgenden Operation von $(h;f) \in H\wr G$ auf den Basiselementen $b_\varphi$
\begin{eqnarray} \label{op}
(h;f).b_\varphi:=f(1).b_{\varphi(h^{-1}(1))}\otimes\ldots \otimes f(n).b_{\varphi(h^{-1}(n))}.
\end{eqnarray}
Die dadurch beschriebene Darstellung von $H\wr G$ wird mit \[\biggl(\stackrel{n}{\textbf{\#}} D\biggr)^\sim\]
bezeichnet, denn die Einschr"ankung der durch (\ref{op}) beschriebenen Darstellung auf $G^*$ ist gerade das $n$-fache "au"sere Tensorprodukt $\textbf{\#}^n D$ von $D$.

{\bf Beispiel 2.3.2} \\
Sei $\left[\lambda\right]=\left[2,2\right]$ Darstellung von $S_4$ und $((132);(12),id,(12)(34)) \in S_3\wr S_4$. 
Weiter sei $\left\{b_1,b_2\right\}$ die $K$-Basis von $V$, so operiert $((132);e)$ wie folgt auf den 8 Basiselementen $\left\{b_{\varphi_1}, \ldots ,b_{\varphi_8}\right\}$ von $\bigotimes^3 V$ :
\[ \begin{array}{ccc}
b_{\varphi_1}:=b_1\otimes b_1\otimes b_1 &\mapsto& b_1\otimes b_1\otimes b_1 = b_{\varphi_1}\\
b_{\varphi_2}:=b_1\otimes b_1\otimes b_2 &\mapsto& b_1\otimes b_2\otimes b_1 = b_{\varphi_3}\\
b_{\varphi_3}:=b_1\otimes b_2\otimes b_1 &\mapsto& b_2\otimes b_1\otimes b_1 = b_{\varphi_5}\\
b_{\varphi_4}:=b_1\otimes b_2\otimes b_2 &\mapsto& b_2\otimes b_2\otimes b_1 = b_{\varphi_7}\\
b_{\varphi_5}:=b_2\otimes b_1\otimes b_1 &\mapsto& b_1\otimes b_1\otimes b_2 = b_{\varphi_2}\\
b_{\varphi_6}:=b_2\otimes b_1\otimes b_2 &\mapsto& b_1\otimes b_2\otimes b_2 = b_{\varphi_4}\\
b_{\varphi_7}:=b_2\otimes b_2\otimes b_1 &\mapsto& b_2\otimes b_1\otimes b_2 = b_{\varphi_6}\\
b_{\varphi_8}:=b_2\otimes b_2\otimes b_2 &\mapsto& b_2\otimes b_2\otimes b_2 = b_{\varphi_8}
\end{array}\]
Die zu dem Element $((132);(12),id,(12)(34))$ unter der Darstellung $({\textbf{\#}^3}\left[2,2\right])^\sim$ korrespondierende Matrix l"asst sich nun folgenderma"sen berechnen:
\begin{eqnarray*}\lefteqn{\biggl(\stackrel{3}{\textbf{\#}}\left[2,2\right]\biggr)^\sim ((132);(12),id,(12)(34))}\\ & &  =\biggl(\stackrel{3}{\textbf{\#}}\left[2,2\right]\biggr)^\sim (1_{S_3};(12),id,(12)(34)) \biggl(\stackrel{3}{\textbf{\#}}\left[2,2\right]\biggr)^\sim ((132);e).
\end{eqnarray*}
Mit der angegebenen Operation von $((132);e)$ auf den Basiselementen $\left\{b_{\varphi_1}, \ldots ,b_{\varphi_8}\right\}$ erh"alt man
\begin{eqnarray*}
\biggl(\stackrel{3}{\textbf{\#}}\left[2,2\right]\biggr)^\sim ((132);e) = \begin{pmatrix} 1 & 0 & 0 & 0 & 0 & 0 & 0 & 0 \\ 0 & 0 & 1 & 0 & 0 & 0 & 0 & 0 \\ 0 & 0 & 0 & 0 & 1 & 0 & 0 & 0 \\ 0 & 0 & 0 & 0 & 0 & 0 & 1 & 0 \\ 0 & 1 & 0 & 0 & 0 & 0 & 0 & 0 \\ 0 & 0 & 0 & 1 & 0 & 0 & 0 & 0 \\ 0 & 0 & 0 & 0 & 0 & 1 & 0 & 0 \\ 0 & 0 & 0 & 0 & 0 & 0 & 0 & 1 \end{pmatrix}
\end{eqnarray*}
Und nach Beispiel 2.2.7 gilt
\begin{eqnarray*}
\biggl(\stackrel{3}{\textbf{\#}}\left[2,2\right]\biggr)^\sim (1_{S_3};(12),id,(12)(34)) = \begin{pmatrix}  1 & 0 \\  0 & -1 \end{pmatrix}\otimes \begin{pmatrix}  1 & 0 \\  0 & 1 \end{pmatrix}\otimes  \begin{pmatrix} 1 & 0 \\  0 & 1 \end{pmatrix} \\ =\begin{pmatrix} 1 & 0 & 0 & 0 & 0 & 0 & 0 & 0 \\ 0 & 1 & 0 & 0 & 0 & 0 & 0 & 0 \\ 0 & 0 & 1 & 0 & 0 & 0 & 0 & 0 \\ 0 & 0 & 0 & 1 & 0 & 0 & 0 & 0 \\ 0 & 0 & 0 & 0 & -1 & 0 & 0 & 0 \\ 0 & 0 & 0 & 0 & 0 & -1 & 0 & 0 \\ 0 & 0 & 0 & 0 & 0 & 0 & -1 & 0 \\ 0 & 0 & 0 & 0 & 0 & 0 & 0 & -1 \end{pmatrix}
\end{eqnarray*}
Daraus folgt:
\begin{eqnarray*}
\biggl(\stackrel{3}{\textbf{\#}}\left[2,2\right]\biggr)^\sim ((132);(12),id,(12)(34)) = \begin{pmatrix} 1 & 0 & 0 & 0 & 0 & 0 & 0 & 0 \\ 0 & 0 & 1 & 0 & 0 & 0 & 0 & 0 \\ 0 & 0 & 0 & 0 & -1 & 0 & 0 & 0 \\ 0 & 0 & 0 & 0 & 0 & 0 & -1 & 0 \\ 0 & 1 & 0 & 0 & 0 & 0 & 0 & 0 \\ 0 & 0 & 0 & 1 & 0 & 0 & 0 & 0 \\ 0 & 0 & 0 & 0 & 0 & -1 & 0 & 0 \\ 0 & 0 & 0 & 0 & 0 & 0 & 0 & -1 \end{pmatrix}
\end{eqnarray*}
Das untenstehende Lemma 2.3.3 aus \cite[S.149]{jk} gibt an, wie man die zur Darstellung $(\textbf{\#}^n D )^\sim$ geh"origen Charaktere berechnet. Dabei sei $g_\nu(h;f)$ die durch Gleichung (\ref{zykelprodukt}) definierte Zykeldarstellung von $(h;f)$. Der Beweis dieses Lemmas wird hier angegeben, da sich die Notation von der in \cite{jk} in einigen Punkten unterscheidet.

{\bf Lemma 2.3.3} \begin{it} 
Es gilt:
\begin{eqnarray} \label{char}
\chi^{(\textbf{\#}^n D)^\sim}(h;f)= \prod_{\nu=1}^{c(h)} \chi^D(g_\nu(h;f)) \ \forall (h;f) \in H\wr G.
\end{eqnarray}\end{it}
\bew
Nach (\ref{VR}) kann die Spur von $(h;f)$ durch
\[\chi^{(\textbf{\#}^n D)^\sim}(h;f)=\sum_\varphi\left[ \mbox{Koeffizient von} \ b_\varphi \  \mbox{in} \  (h;f).b_\varphi\right]\]
berechnet werden. Sei $\mathbb{D}$ die zur Darstellung $D$ korrespondierende Matrix mit \\ $\mathbb{D}(g)=:(d_{ik}(g))$, so l"asst sich der Koeffizient mit Gleichung (\ref{op}) wie folgt berechnen:
\begin{eqnarray*}
(h;f).b_\varphi &=& \Bigl(\sum_{i_1=1}^m d_{i_1\varphi(h^{-1}(1))}(f(1))b_{i_1}\Bigr)\otimes \cdots \otimes \Bigl(\sum_{i_n=1}^m d_{i_n\varphi(h^{-1}(n))}(f(n))b_{i_n}\Bigr)\\
&=& \sum_{1\leq i_1,\ldots ,i_n\leq m} d_{i_1\varphi(h^{-1}(1))}(f(1)) \cdots d_{i_n\varphi(h^{-1}(n))}(f(n))b_{i_1} \otimes \cdots \otimes b_{i_n}\\
&=& \sum_\psi\Bigl( \prod_j d_{\psi(j)\varphi(h^{-1}(j))}(f(j))\Bigr)b_\psi.
\end{eqnarray*}
Der Koeffizient von $b_\varphi$ ist also $\prod_{j=1}^n d_{\varphi(j)\varphi(h^{-1}(j))}(f(j))$. Sortiert man nun die Faktoren in geeigneter Weise um, so erh"alt man das Produkt
\begin{eqnarray*}
\prod_{\nu=1}^{c(h)} d_{\varphi(j_\nu)\varphi(h^{-1}(j_\nu))}(f(j_\nu))d_{\varphi(h^{-1}(j_\nu))\varphi(h^{-2}(j_\nu))}(f(h^{-1}(j_\nu)))\cdots d_{\varphi(h^{l_\nu-1}(j_\nu))\varphi(j_\nu)}(f(h^{l_\nu-1}(j_\nu))).
\end{eqnarray*}
Summiert man diesen Ausdruck noch "uber alle Abbildungen $\varphi \in [m]^{[n]}$, so f"uhrt dies zur gew"unschten Charaktergleichung.\qed 

Im Weiteren soll nun erl"autert werden, wie man ein vollst"andiges System paarweiser nicht "aquivalenter irreduzibler Darstellungen von $H\wr G$ erh"alt. Dazu betrachte man ein vollst"andiges System $D^1,\ldots ,D^r$ paarweiser nicht "aquivalenter irreduzibler Darstellungen von $G$ "uber $K$ mit zugeh"origen Vektorr"aumen $V^j$. Die irreduziblen Darstellungen von $G^*$ sind dann von der Form $D^*:=D_1${\footnotesize{\textbf \#}} $\cdots$ {\footnotesize{\textbf \#}} $D_n$ ,wobei $D_i \in \left\{D^1,\ldots ,D^r\right\}$. Der zugrundeliegende Vektorraum ist $\bigotimes\limits_i V := V_1\otimes_K \ldots \otimes_K V_n$, wobei $V_i := V^j$, falls $D_i=D^j$. 

{\bf Definition 2.3.4} \\
a) Sei $n_j$ die Anzahl der Faktoren $D_i$ von $D^*$, welche mit $D^j$ "ubereinstimmen, $1\leq j\leq r$, \hspace*{0.4cm} so hei"st $(n):=(n_1,\ldots ,n_r)$ {\it Typ von $D^*$}.\\
b) Die {\it Tr"agheitsgruppe von $D^*$} ist definiert als \[H_{D^*}\wr G := \left\{(h;f)| D^{*(h;f)}\sim D^*\right\},\]
\hspace*{0.4cm} wobei $\sim$ die "Aquivalenz der Darstellungen bedeutet und $D^{*(h;f)}$ die zu $D^*$ konjugierte \hspace*{0.4cm} Darstellung ist, das hei"st 
\[D^{*(h;f)}(1_H;f')= D^*(h;f)^{-1}(1_H;f')(h;f).\]
{\bf Lemma 2.3.5} \\
\begin{it}
Sei $(n)$ der Typ von $D^*$ und sei $S_{(n)}:= S_{n_1}\times \ldots \times S_{n_r}$, so gilt f"ur die Tr"agheitsgruppe 
\[H_{D^*}\wr G= H\cap S_{(n)}\wr G.\] \end{it}
\bew siehe \cite[S.152, Lemma 4.3.27]{jk}.\qed 

Nach Cliffords Theorie lassen sich die irreduziblen Darstellungen $D^*$ zu Darstellungen von $H_{D^*}\wr G$ fortsetzen. Diese fortgesetzten Darstellungen werden mit $\widetilde{D^*}$ bezeichnet. \\
Sei $D''$ eine irreduzible Darstellung von $H_{D^*}$, so erh"alt man eine zweite irreduzible Darstellung $D'$ von $H_{D^*}\wr G$ durch 
\[ D'(h;f):= D''(h).\]
Die Multiplikation dieser beiden Darstellungen ergibt eine dritte irreduzible Darstellung von $H_{D^*}\wr G$: Das innere Tensorprodukt 
\[D'\otimes \widetilde{D^*}.\]
Cliffords Theorie "uber die Darstellungen von Gruppen mit normalen Untergruppen liefert das Resultat, dass jede irreduzible Darstellung von $H\wr G$ von der Form \[D:= (D'\otimes \widetilde{D^*}) \uparrow H\wr G\] ist und besagt im folgenden Theorem zudem, wie man ein vollst"andiges System irreduzibler Darstellungen von $H\wr G$ erh"alt. 

{\bf Theorem 2.3.6} \\ \begin{it}
Die irreduzible K-Darstellung $D:=(D'\otimes \widetilde{D^*})\uparrow H\wr G$ durchl"auft ein vollst"andiges System paarweise nicht "aquivalenter irreduzibler Darstellungen von $H\wr G$, falls $D^*$ ein vollst"andiges System paarweise nicht konjugierter (bez"uglich $H'$) aber irreduzibler K-Darstellungen von $G^*$ und 
$D''$ bei festem $D^*$ ein vollst"andiges System paarweise nicht "aquivalenter und irreduzibler K-Darsellungen von $H\cap S_{(n)}$ durchl"auft. \end{it}

{\bf Beispiel 2.3.7} \\
Die irreduziblen Darstellungen von $S_3$ sind $[3], [2,1]$ und $[1^3]$, die von $S_2$ sind $[2]$ und $[1^2]$ und die von $S_1$ ist $[1]$. Ein vollst"andiges System irreduzibler Darstellungen von $S_3\wr S_3$ ist dann nach Theorem 2.3.6
\[
\begin{array}{lll}
([3];[3]) & ([2,1];[3]) & ([1^3];[3]) \\
([3];[2,1]) & ([2,1];[2,1]) & ([1^3];[2,1]) \\
([3];[1^3]) & ([2,1];[1^3]) & ([1^3];[1^3]) \\
\end{array} \]
\[
\begin{array}{ll}
([2]'\otimes ([3]${\footnotesize\textbf{\#}}$[3]${\footnotesize\textbf{\#}}$[2,1]))\uparrow S_3\wr S_3 & ([1^2]'\otimes ([3]${\footnotesize\textbf{\#}}$[3]${\footnotesize\textbf{\#}}$[2,1]))\uparrow S_3\wr S_3  \\
([2]'\otimes ([3]${\footnotesize\textbf{\#}}$[2,1]${\footnotesize\textbf{\#}}$[2,1]))\uparrow S_3\wr S_3 & ([1^2]'\otimes ([3]${\footnotesize\textbf{\#}}$[2,1]${\footnotesize\textbf{\#}}$[2,1]))\uparrow S_3\wr S_3 \\
([2]'\otimes ([3]${\footnotesize\textbf{\#}}$[3]${\footnotesize\textbf{\#}}$[1^3]))\uparrow S_3\wr S_3 & ([1^2]'\otimes ([3]${\footnotesize\textbf{\#}}$[3]${\footnotesize\textbf{\#}}$[1^3]))\uparrow S_3\wr S_3  \\
([2]'\otimes ([3]${\footnotesize\textbf{\#}}$[1^3]${\footnotesize\textbf{\#}}$[1^3]))\uparrow S_3\wr S_3 & ([1^2]'\otimes ([3]${\footnotesize\textbf{\#}}$[1^3]${\footnotesize\textbf{\#}}$[1^3]))\uparrow S_3\wr S_3 \\
([2]'\otimes ([1^3]${\footnotesize\textbf{\#}}$[2,1]${\footnotesize\textbf{\#}}$[2,1]))\uparrow S_3\wr S_3 & ([1^2]'\otimes ([1^3]${\footnotesize\textbf{\#}}$[2,1]${\footnotesize\textbf{\#}}$[2,1]))\uparrow S_3\wr S_3  \\
([2]'\otimes ([2,1]${\footnotesize\textbf{\#}}$[1^3]${\footnotesize\textbf{\#}}$[1^3]))\uparrow S_3\wr S_3 & ([1^2]'\otimes ([2,1]${\footnotesize\textbf{\#}}$[1^3]${\footnotesize\textbf{\#}}$[1^3]))\uparrow S_3\wr S_3  \\
\end{array}
\]
\begin{center}
$([1]'\otimes ([3]${\footnotesize\textbf{\#}}$[2,1]${\footnotesize\textbf{\#}}$[1^3]))\uparrow S_3\wr S_3$
\end{center}
wobei $([\lambda];[\mu]):= ([\lambda]\otimes (\stackrel{3}{\textbf{\#}}[\mu])^\sim)$. 

{\bf Lemma 2.3.8} \begin{it} 
Es gilt:
\begin{eqnarray} \label{char2}
\chi^{\widetilde{D^*}}(h;f)= \prod_{\nu=1}^{c(h)} \chi^{D_{j_\nu}}(g_\nu(h;f)) \ \forall (h;f) \in H_{D^*}\wr G.
\end{eqnarray}\end{it}
\bew
Analog zum Beweis von Lemma 2.3.3 erh"alt man anstelle der letzten Gleichung im Beweis jetzt die Gleichung
\begin{eqnarray*}
\prod_{\nu=1}^{c(h)} d^{j_\nu}_{\varphi(j_\nu)\varphi(h^{-1}(j_\nu))}(f(j_\nu))d^{h^{-1}(j_\nu)}_{\varphi(h^{-1}(j_\nu))\varphi(h^{-2}(j_\nu))}(f(h^{-1}(j_\nu)))\cdots d^{h^{l_\nu-1}(j_\nu)}_{\varphi(h^{l_\nu-1}(j_\nu))\varphi(j_\nu)}(f(h^{l_\nu-1}(j_\nu))).
\end{eqnarray*}
Da $h \in H_{D^*}$ gilt $d^{j_\nu}=d^{h^{-1}(j_\nu)}= \cdots = d^{h^{l_\nu-1}(j_\nu)}$ und man erh"alt wieder durch Summation "uber die Abbildungen $\varphi$ die gew"unschte Charaktergleichung.\qed

\textsc{Bemerkung:}\\
Die Abbildungen $\varphi$ in Lemma 2.3.8 stimmen nicht mit denen in Lemma 2.3.3 "uberein. Da die zugrundeliegenden Vektorr"aume der Darstellungen $D_i$ nicht identisch sein m"ussen.

{\bf Beispiel 2.3.9} \\
Betrachte die Darstellung $\widetilde{D^*}=([2,1]${\footnotesize{\textbf \#}}$[2,1]${\footnotesize{\textbf \#}}$[1^3])$ und das Element $((12);(12),id,(123)) \in S_{3_{D^*}}\wr S_3$, so gilt:
\begin{eqnarray*}
[2,1]\#[2,1]\#[1^3](1;(12),id,(123)) & = & \begin{pmatrix}  -1 & 0 \\  0 & 1 \end{pmatrix}\otimes \begin{pmatrix}  1 & 0 \\  0 & 1 \end{pmatrix}\otimes  \begin{pmatrix} -1  \end{pmatrix} \\ & = & \begin{pmatrix} 1 & 0 & 0 & 0 \\ 0 & 1 & 0 & 0 \\ 0 & 0 & -1 & 0 \\ 0 & 0 & 0 & -1 \end{pmatrix}.
\end{eqnarray*}
Sei $\left\{b_1,b_2\right\}$ Basis von $V_1=V_2$ und $\left\{c_1\right\}$ Basis von $V_3$, so operiert $((12),e)$ wie folgt auf den 4 Basiselementen $\left\{b_{\varphi_1},\ldots ,b_{\varphi_4}\right\}$
\[ \begin{array}{ccc}
b_{\varphi_1}:=b_1\otimes b_1\otimes c_1 &\mapsto& b_1\otimes b_1\otimes c_1 = b_{\varphi_1}\\
b_{\varphi_2}:=b_1\otimes b_2\otimes c_1 &\mapsto& b_2\otimes b_1\otimes c_1 = b_{\varphi_3}\\
b_{\varphi_3}:=b_2\otimes b_1\otimes c_1 &\mapsto& b_1\otimes b_2\otimes c_1 = b_{\varphi_2}\\
b_{\varphi_4}:=b_2\otimes b_2\otimes c_1 &\mapsto& b_2\otimes b_2\otimes c_1 = b_{\varphi_4}
\end{array}\]
Somit folgt
\begin{eqnarray*}
[2,1]\#[2,1]\#[1^3]((12);(12),id,(123)) = \begin{pmatrix} 1 & 0 & 0 & 0 \\ 0 & 0 & 1 & 0 \\ 0 & -1 & 0 & 0 \\ 0 & 0 & 0 & -1 \end{pmatrix}
\end{eqnarray*}
und $\chi^{\widetilde{D^*}}((12);(12),id,(123)) = 0$.

\textsc{Charakterberechnung f"ur induzierte Charaktere:}\\
Die Formel f"ur die Berechnung der induzierten Charakterwerte (vgl. Formel (\ref{ind})) ist gegeben durch
\begin{eqnarray}\label{1}
\chi^D(h;f)=\begin{cases} 
\frac{1}{|H_{D^*}\wr G|} \sum\limits_{x\in H\wr G} \dot{\chi}^{D'\otimes \widetilde{D^*}}(x(h;f)x^{-1}),&\mbox{falls} \ H_{D^*} \neq H \\
\chi^{D'\otimes \widetilde{D^*}}(h;f), & \mbox{sonst} 
\end{cases} \end{eqnarray}
\begin{eqnarray} \label{2} \mbox{wobei} \ \dot{\chi}^{D'\otimes \widetilde{D^*}}(h';f')=\begin{cases} \chi^{D'\otimes \widetilde{D^*}}(h';f'), & \mbox{falls} \ (h';f') \in H_{D^*}\wr G \\
0, & \mbox{sonst} \end{cases} \end{eqnarray}
\bq \label{3}\mbox{und} \ \chi^{D'\otimes \widetilde{D^*}}(h;f)&=&\chi^{D'}(h;f)\chi^{\widetilde{D^*}}(h;f) \nonumber\\
&=&\chi^{D''}(h)\chi^{\widetilde{D^*}}(h;f)\nonumber\\
&=&\chi^{D''}(h)\prod_{\nu=1}^{c(h)} \chi^{D_{j_\nu}}(g_\nu(h;f))  
\eq
Dabei ist $D_{j_\nu}=D$ f"ur den Fall, dass $D^*=\#^n D$ ist.
\chapter{\textbf{ $C_4\wr S_n$}}
Dieses Kapitel soll zum einen die Betrachtung des Kranzproduktes $C_4\wr S_n$ f"ur die Berechnung des Geschlechts eines Origamis $O$ motivieren. Zum anderen werden die Resultate der Charaktertheorie aus Kapitel 2 verwendet, um die Verteilung des Geschlechts eines zuf"allig gew"ahlten Origami, das hei"st eines Origami mit zuf"alligen Verklebungen, zu berechnen.
\section{Motivation}
Die Berechung des Geschlechts $g$ eines zuf"allig gew"ahlten Origami, bestehend aus $n$-Kopien des euklidischen Einheitsquadrates, geht nach der Euler-Charakteristik $\chi= E-K+F$ und der Beziehung $\chi=2-2g$ zur"uck auf die Berechnung der Eckenanzahl nach Verklebung. Die Anzahl der Fl"achen ist $n$ und die Anzahl der Kanten ist $2n$.

Betrachtet man die zwei Elemente
\[ \sigma:=((13)(24);\sigma_a,\sigma_b,\sigma_a^{-1},\sigma_b^{-1}) \ \mbox{und} \ \tau:=((1234); id,id,id,id) \in C_4\wr S_n,\] 
so stimmt die Anzahl der Identifikationsklassen der Ecken nach der Verklebung mit der Anzahl der Bahnen von folgender Operation "uberein:
\begin{eqnarray} \label{opbahn}
\phi: \left\langle \sigma\tau\right\rangle \times ([4]\times[n]) \rightarrow ([4]\times[n])\nonumber \\
\phi((\sigma\tau)^k,(i,j))=(\tilde{\pi}^k(i),(f_k)_{\tilde{\pi}^k(i)}(j))
\end{eqnarray}
wobei $\left\langle \sigma\tau\right\rangle \ni (\sigma\tau)^k=(\tilde{\pi}^k;(f_k)_1,(f_k)_2,(f_k)_3,(f_k)_4)$, $\tilde{\pi}=(1432)$, $(f_k)_i \in S_n$ und $[n]:=\left\{1,\ldots ,n\right\}$ ist. \\ Hierbei entspricht $(i,j) \in [4]\times[n]$ der $i$-ten Kante im $j$-ten euklidischen Einheitsquadrat, so dass die zu verklebenden Einheitsquadrate zuk"unftig immer die untenstehende Kantenbeschriftung haben, auch wenn diese nicht explizit angegeben ist. Ebenso wird ab jetzt anstelle von $\phi((\sigma\tau)^k,(i,j))$ die Notation $(\sigma\tau)^k.(i,j)$ verwendet.
\begin{center}
\begin{picture}(7.5,1)
\multiput(0,0)(2,0){2}{\line (1,0){1.5}}
\multiput(0,0)(1.5,0){2}{\line (0,-1){1.5}}
\multiput(2,0)(1.5,0){2}{\line (0,-1){1.5}}
\multiput(0,-1.5)(2,0){2}{\line (1,0){1.5}}
\multiput(6,0)(0,-1.5){2}{\line (1,0){1.5}}
\multiput(6,0)(1.5,0){2}{\line (0,-1){1.5}}
\put(0.5,-0.8){ $1$}
\put(2.5,-0.8){ $2$}
\put(6.5,-0.8){ $n$}
\multiput(4,-0.8)(0.7,0){3}{\circle*{0.08}}
\multiput(0.65,-0.2)(2,0){2}{\scriptsize $4$}
\put(6.65,-0.2){\scriptsize $4$}
\multiput(0.08,-0.8)(2,0){2}{\scriptsize $1$}
\put(6.08,-0.8){\scriptsize $1$}
\multiput(1.34,-0.8)(2,0){2}{\scriptsize $3$}
\put(7.34,-0.8){\scriptsize $3$}
\multiput(0.65,-1.44)(2,0){2}{\scriptsize $2$}
\put(6.65,-1.44){\scriptsize $2$}
\end{picture}
\end{center}
\vspace{0.5cm}
\begin{center}
{\small Abb.4: $n$ Euklidische Einheitsquadrate mit Kantenbeschriftung.}
\end{center}
Die Beispiele 3.1.1 a) und b) werden zeigen, dass die Betrachtung des Kranzproduktes im Falle von Origamis, die durch Verklebungsvorschriften definiert werden, sehr anschaulich ist, denn bildlich bedeutet $\tau$ eine Vierteldrehung nach rechts und $\sigma$ ein Wechsel zur verklebten Kante. Zeichnet man nun alle Drehungen ein, so umrundet man alle Ecken der selben Markierung, indem man alle Kanten einer Bahn durchl"auft.

{\bf Beispiel 3.1.1}
\begin{enumerate}
\item[a)]
Bei dem Origami $O_1$ aus Beispiel 1.1.2 b) ist $\sigma=((13)(24);(12),(12),(12),(12))$ und $\tau=((1234);id,id,id,id)$, so dass $\sigma\tau=((1432);(12),(12),(12),(12))$. Man erh"alt unter obiger Operateration zwei Bahnen
\begin{eqnarray*}
\sigma\tau.(1,1)=(4,2),\ \sigma\tau.(4,2)=(3,1),\ \sigma\tau.(3,1)=(2,2),\ \sigma\tau.(2,2)=(1,1)\\
\sigma\tau.(1,2)=(4,1),\ \sigma\tau.(4,1)=(3,2),\ \sigma\tau.(3,2)=(2,1),\ \sigma\tau.(2,1)=(1,2).
\end{eqnarray*} 
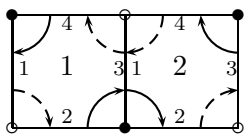
\begin{center}
\begin{pspicture}(2,-1)
\multiput(0,0)(1.5,0){2}{\line (1,0){1.5}}
\multiput(0,0)(1.5,0){2}{\line (0,-1){1.5}}
\multiput(1.5,0)(1.5,0){2}{\line (0,-1){1.5}}
\multiput(0,-1.5)(1.5,0){2}{\line (1,0){1.5}}
\put(0.5,-0.8){ $1$}
\put(2,-0.8){ $2$}
\multiput(0.65,-0.2)(1.5,0){2}{\scriptsize $4$}
\multiput(0.08,-0.8)(1.5,0){2}{\scriptsize $1$}
\multiput(1.34,-0.8)(1.5,0){2}{\scriptsize $3$}
\multiput(0.65,-1.44)(1.5,0){2}{\scriptsize $2$}
\multiput(-0.1,-0.1)(3,0){2}{$\bullet$}
\put(1.4,-0.1){$\circ$}
\multiput(-0.1,-1.6)(3,0){2}{$\circ$}
\put(1.4,-1.6){$\bullet$}
\put(0,-1.5){\psarc[linestyle=dashed]{<-}{0.5}{0}{90}}
\put(3,-1.5){\psarc[linestyle=dashed]{<-}{0.5}{90}{180}}
\put(1.5,0){\psarc[linestyle=dashed]{<-}{0.5}{180}{270}}
\put(1.5,0){\psarc[linestyle=dashed]{<-}{0.5}{270}{0}}
\put(0,0){\psarc{<-}{0.5}{270}{0}}
\put(3,0){\psarc{<-}{0.5}{180}{270}}
\put(1.5,-1.5){\psarc{<-}{0.5}{90}{180}}
\put(1.5,-1.5){\psarc{<-}{0.5}{0}{90}}
\end{pspicture}
\end{center}
\vspace{0.5cm}
\item[b)]
F"ur das Origami $D$ erh"alt man mit $\sigma_a=(123)$ und $\sigma_b=(145)(23)$ die folgenden drei Bahnen
\begin{eqnarray*}\left\langle \sigma\tau\right\rangle.(1,1)&=&\left\{(4,5),(3,5),(2,1),(1,2),(4,3),(3,2),(2,3),(1,1)\right\} \\
\left\langle \sigma\tau\right\rangle.(1,3)&=&\left\{(4,2),(3,1),(2,4),(1,4),(4,1),(3,3),(2,2),(1,3)\right\} \\
\left\langle \sigma\tau\right\rangle.(1,5)&=&\left\{(4,4),(3,4),(2,5),(1,5)\right\}.
\end{eqnarray*}
\end{enumerate}
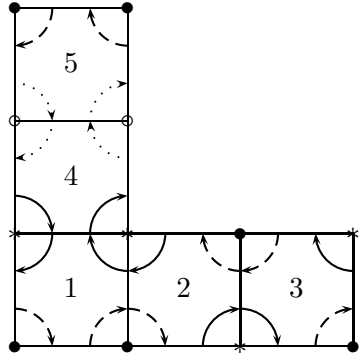
\begin{center}
\begin{picture}(3,3)
\multiput(0,0)(1.5,0){3}{\line (1,0){1.5}}
\multiput(0,0)(1.5,0){4}{\line (0,-1){1.5}}
\multiput(0,-1.5)(1.5,0){3}{\line (1,0){1.5}}
\multiput(0,0)(1.5,0){2}{\line (0,1){1.5}}
\multiput(0,1.5)(1.5,0){2}{\line (0,1){1.5}}
\multiput(0,1.5)(0,1.5){2}{\line (1,0){1.5}}
\put(0.65,-0.85){$1$}
\put(2.15,-0.85){$2$}
\put(3.65,-0.85){$3$}
\put(0.65,0.65){$4$}
\put(0.65,2.15){$5$}
\multiput(-0.1,-0.1)(1.5,0){2}{$\ast$}
\put(2.9,-0.1){$\bullet$}
\put(4.4,-0.1){$\ast$}
\multiput(-0.1,-1.6)(1.5,0){2}{$\bullet$}
\put(2.9,-1.6){$\ast$}
\put(4.4,-1.6){$\bullet$}
\multiput(-0.1,1.4)(1.5,0){2}{$\circ$}
\multiput(-0.1,2.9)(1.5,0){2}{$\bullet$}
\put(0,-1.5){\psarc[linestyle=dashed]{<-}{0.5}{0}{90}}
\put(1.5,-1.5){\psarc[linestyle=dashed]{<-}{0.5}{0}{90}}
\put(1.5,-1.5){\psarc[linestyle=dashed]{<-}{0.5}{90}{180}}
\put(3,0){\psarc[linestyle=dashed]{<-}{0.5}{180}{270}}
\put(3,0){\psarc[linestyle=dashed]{<-}{0.5}{270}{00}}
\put(4.5,-1.5){\psarc[linestyle=dashed]{<-}{0.5}{90}{180}}
\put(0,3){\psarc[linestyle=dashed]{<-}{0.5}{270}{0}}
\put(1.5,3){\psarc[linestyle=dashed]{<-}{0.5}{180}{270}}
\put(0,0){\psarc{<-}{0.5}{270}{0}}
\put(0,0){\psarc{<-}{0.5}{0}{90}}
\put(1.5,0){\psarc{<-}{0.5}{90}{180}}
\put(1.5,0){\psarc{<-}{0.5}{180}{270}}
\put(1.5,0){\psarc{<-}{0.5}{270}{0}}
\put(3,-1.5){\psarc{<-}{0.5}{0}{90}}
\put(3,-1.5){\psarc{<-}{0.5}{90}{180}}
\put(4.5,0){\psarc{<-}{0.5}{180}{270}}
\put(0,1.5){\psarc[linestyle=dotted]{<-}{0.5}{0}{90}}
\put(1.5,1.5){\psarc[linestyle=dotted]{<-}{0.5}{90}{180}}
\put(1.5,1.5){\psarc[linestyle=dotted]{<-}{0.5}{180}{270}}
\put(0,1.5){\psarc[linestyle=dotted]{<-}{0.5}{270}{0}}
\end{picture}
\end{center}
\vspace{2cm}
Eine weitere M"oglichkeit, die Anzahl der Ecken nach Verklebung zu berechnen, ist die Betrachtung der Anzahl der Zykel des Elements $\sigma_a\sigma_b\sigma_a^{-1}\sigma_b^{-1} \in S_n$. Dies h"atte zur Folge, dass nicht Formel (\ref{frob}), sondern die Anzahl der L"osungen $N^{S_n}(z)$ der Gleichung $[x,y]=x^{-1}y^{-1}xy=z$ betrachtet werden m"usste. Diese Anzahl ist gegeben durch
\[ N^{S_n}(z)= n!\sum_\chi \frac{\chi(z)}{\chi(1)}\]
(vgl \cite[S.625]{pm2002}). Von diesem Ansatz soll hier jedoch nicht ausgegangen werden, da es auch Ziel dieser Arbeit ist, das Kranzprodukt und dessen Darstellung genauer zu untersuchen.
\section{Wahrscheinlichkeit der Eckenanzahl}
Nach 3.1 l"asst sich die Wahrscheinlichkeit der Eckenanzahl nach Identifikation "uber die Anzahl der Bahnen, der durch (\ref{opbahn}) definierten Operation, berechnen, das hei"st
\[ P(\# \mbox{Ecken nach Identifiktaion} =k)=P(\# \mbox{Bahnen} =k).\]
{\bf Behauptung 3.2.1} \begin{it} 
\[P(\# \mbox{Bahnen} =k)=\sum_{\pi:\pi' \ hat \ k \ Zykel} P(\sigma\tau=\pi),\]
wobei $\pi=(\tilde{\pi};\pi_1,\pi_2,\pi_3,\pi_4) \in C_u\wr S_n$ und $\pi'=\pi_1\circ\pi_2\circ\pi_3\circ\pi_4$ ist.\end{it} \\
\bew
Mit $\pi=(\tilde{\pi};\pi_1,\pi_2,\pi_3,\pi_4)$ und $\tilde\pi=(1432)$ folgt f"ur die Operation von $\pi$ auf $(i,j)$, dass
\[\pi.(i,j)=((1432)(i);\pi_{(1432)(i)}(j)).\]
Dies impliziert, dass die Bahnen dieser Operation immer 4$\cdot$b, b $\in \mathbb{N}$, Elemente haben, da die Ordnung von $\tilde{\pi}$ vier ist.\\ Betrachtet man nun \OE \ das Tupel $(1,j)$, so ist $\pi^4(1,j)=(1;\pi_1\circ\pi_2\circ\pi_3\circ\pi_4(j))$. Liegt $j$ nun in einem Zykel der L"ange b von $\tilde{\pi}$, so hat die Bahn 4$\cdot$b Elemente. Die Anzahl der Bahnen enspricht dann genau der Anzahl der Zykel von $\tilde{\pi}$.\\
(W"ahlt man ein anderes Tupel, so erh"alt man nicht $\tilde{\pi}$, sondern ein zu $\tilde{\pi}$ konjugiertes Element in $S_n$.) \qed

{\bf Theorem 3.2.2} \begin{it} \\
Sei die Anzahl $m$ der Zykel von $\pi' \leq \sqrt[3]{n}$, so gilt
\begin{eqnarray*} \label{Thm 3.2.2} P(\sigma\tau=\pi)=\frac{1}{n!^4}(1+sgn(\pi')+\mathcal{O}(n^{-\frac{2}{3}})). \end{eqnarray*}
\end{it}\\
{\bf Korallar 3.2.3}\begin{it}\\ 
Es gilt
\bq \label{Kor 3.2.3} P(\sigma\tau=\pi)= \frac{1}{|C_4\wr S_n|}\sum_D \frac{\chi^D(\sigma)\chi^D(\tau)\chi^D(\pi)}{\chi^D(1)},\eq
wobei die Summe "uber alle irreduziblen Darstellungen  D von $C_4\wr S_n$ l"auft.\end{it}\\
\bew 
Nach (\ref{frob}) gilt, dass die Anzahl der L"osungen der Gleichung $\sigma\tau=\pi$ gleich $\frac{|C_\sigma||C_\tau|}{|G|} \sum_{\chi \in Irr(G)} \frac{\chi(\sigma)\chi(\tau)\chi(\pi^{-1})}{\chi(1)}$ ist, wobei $C_\sigma$ bzw. $C_\tau$ die zu $\sigma$ bzw. $\tau$ geh"origen Konjugationsklassen sind. Da nach Lemma 2.1.10 $(iii)$ der Charakter eine Klassenfunktion ist, kann man anstelle von $\chi(C_\sigma)$ bzw. $\chi(C_\tau)$ auch $\chi(\sigma)$ bzw.$\chi(\tau)$ schreiben. Teilt man noch durch die Anzahl der m"oglichen Permutationen von $\sigma\tau$, so erh"alt man die Aussage aus dem Korollar.\qed 

Zun"achst sollten also die Charaktere der Element $\sigma$ und $\tau$ f"ur den Beweis des Theorems 3.2.2 bestimmt werden. Hierf"ur betrachte man die Gleichungen (\ref{1}),(\ref{2}) und (\ref{3}).\\
Aus (\ref{2}) folgt, dass falls $\dot{\chi}(\tau) \neq 0$ sein soll, so muss gelten $D^*=\lambda${\footnotesize{\textbf \#}}$\lambda${\footnotesize{\textbf \#}}$\lambda${\footnotesize{\textbf \#}}$\lambda$, wobei $\lambda$ eine irreduzible Darstellung von $S_n$ ist. Es folgt also, dass ${C_4}_{D^*}=C_4\cap S_4=C_4$ ist.\\ \\
Die Charaktere lassen sich demnach durch Gleichung (\ref{3}) mit dem Zykelprodukt aus Definition 1.2.4 wie folgt darstellen:
\bq \chi^D(\sigma)& = & \chi^{D'\otimes \widetilde{D^*}}(\sigma)=\chi^{D''}((13)(24))\chi^\lambda(1)^2 \\
\chi^D(\tau)& = & \chi^{D''}((1234))\chi^\lambda(1) \eq \newpage
\bq \chi^D(\pi)& = & \chi^{D''}((1234))\chi^\lambda(g_\nu(\pi^{-1})) \nonumber \\
& = & \chi^{D''}((1234))\chi^\lambda(\pi'^{-1}) \nonumber \\
& = & \chi^{D''}((1234))\chi^\lambda(\pi') \ \ \ \ \mbox{nach 2.1.10} \ (ii)\\
\chi^D(1)& = & \chi^{D''}(1)\chi^\lambda(1)^4
\eq

Die Charakterwerte der Elemente von $C_4$ sind durch die untenstehende Tabelle gegeben.
\begin{center}
\begin{tabular}{r||cccc} \hline
    & $id$ & $(1234)$ & $(13)(24)$ & $(1432)$ \\ \hline\hline 
$\chi^1$ & 1 & 1 & 1 & 1 \\
$\chi^2$ & 1 & i & -1 & -i \\
$\chi^3$ & 1 & -1 & 1 & -1 \\
$\chi^4$ & 1 & -i & -1 & i \\ \hline
\end{tabular}
\end{center}

{\bf Korollar 3.2.4}\begin{it}\\
Es gilt
\[ P(\sigma\tau=\pi)= \frac{1}{n!^4}\sum_\lambda \frac{\chi^\lambda(\pi')}{\chi^\lambda(1)},\]
wobei die Summe "uber alle irreduziblen Darstellungen $\lambda$ von $S_n$ l"auft.\end{it} \\
\bew Nach den Gleichungen (3.4)-(3.7) l"asst sich Gleichung (\ref{Kor 3.2.3}) umschreiben zu
\[ P(\sigma\tau=\pi)= \frac{1}{4n!^4}\sum_D \frac{\chi^{D''}((13)(24))\chi^\lambda(1)^2\chi^{D''}((1234))\chi^\lambda(1)\chi^{D''}((1234))\chi^\lambda(\pi')}{\chi^{D''}(1)\chi^\lambda(1)^4} \]
Nach obenstehender Charakterwerttabelle gilt jedoch  \[ \frac{\chi^{D''}((13)(24))\chi^{D''}((1234))^2}{\chi^{D''}(1)}=1 \ \mbox{f"ur alle Darstellungen} \ D'' \ \mbox{von} \ C_4.\]
K"urzt man die 4 Darstellungen von $C_4$, so erh"alt man die Aussage des Korollars.\qed 

Extrahiert man nun die Darstellungen $[n],[1^n],[n-1,1]$ und $[2,1^{n-2}]$ und wendet auf diese die Murnaghan-Nakayama Formel (\ref{Mu-Na}) an, so ergibt sich die Wahrscheinlichkeit 
\begin{eqnarray} \label{summand}
 P(\sigma\tau=\pi)= \frac{1}{n!^4}((1 + sgn(\pi'))(1 + \frac{l-1}{n-1})\sum_\lambda \frac{\chi^\lambda(\pi')}{\chi^\lambda(1)},
\end{eqnarray}
wobei die Summe "uber alle irreduziblen Darstellungen $\lambda$ mit $\lambda \neq [n], [1^n], [n-1,1], [2,1^{n-2}]$ von $S_n$ l"auft und $l$ die Anzahl der Fixpunkte von $\pi'$ ist.

Der Rest dieses Unterkapitels besch"aftigt sich mit dem Beweis des Theorems 3.2.2, der sich in zwei Teile gliedert. In Teil 1 soll f"ur $m\leq \sqrt[3]{n}$ gezeigt werden, dass $\chi^\lambda(\pi') \leq \chi^\lambda(1)^{\frac{1}{2}}$ ist. Mit dieser Aussage l"asst sich dann in Teil 2 der Summationsterm aus Korollar 3.2.4 weiter absch"atzen. Sei daf"ur im folgenden $A=n-\lambda_1$ und $\lambda_1 \geq \|\lambda\|$ und die Partition $(\mu)=(\lambda_2,\ldots, \lambda_h)$.

{\bf Beweis von Theorem 3.2.2:}\\
\textsc{Teil 1:}

{\bf Lemma 3.2.5}\begin{it}\\
Sei $\pi \in S_n$ und sei $m$ die Anzahl der Zykel von $\pi'$. Dann gilt f"ur jede irreduzible Darstellung $\lambda$ von $S_n$
\bq |\chi^\lambda(\pi')| \leq (2n)^{\frac{m}{2}}.\eq \\ \end{it}
\bew siehe \cite[S.568,Prop 2.12]{ls2004}.\qed 

Es gen"ugt also $\chi^{\lambda}(1)^{\frac{1}{2}}\geq (2n)^{\frac{\sqrt[3]{n}}{2}}$ zu zeigen.

\fbox{{\bf Fall 1: $\lambda_1 \leq \frac{3}{4}n$}} \\ 
Mit Lemma 8 aus \cite{pm2007} folgt $\chi^\lambda(1) \geq 2^{\frac{n}{8}}$.
Und f"ur $n$ gro"s genug (ca. ab $n$=790) gilt $2^{\frac{n}{16}} \geq (2n)^{\frac{\sqrt[3]{n}}{2}}$.

\fbox{{\bf Fall 2: $\lambda_1 \geq \frac{3}{4}n$ und $A>\sqrt{n}$ }} \\ 
F"ur diesen Fall betrachte man das zu der Partition $(\lambda)=(\lambda_1, \ldots ,\lambda_h)$ geh"orige Diagramm 
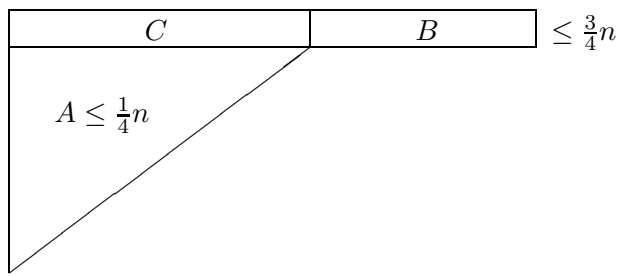
\begin{center} \begin{picture}(7,-7)
\put(0,0){\line(1,0){7}}
\put(0,0){\line(0,-1){3.5}}
\put(0,-0.5){\line(1,0){7}}
\put(4,0){\line(0,-1){0.5}}
\put(7,0){\line(0,-1){0.5}}
\put(0,-3.5){\line(4,3){4}}
\put(1.8,-0.4){$C$}
\put(5.4,-0.4){$B$}
\put(0.6,-1.5){$A\leq \frac{1}{4}n$}
\put(7.2,-0.4){$\leq \frac{3}{4}n$}
\end{picture}
\end{center}
\vspace{3 cm}
wobei $A$ die Zeilen $\lambda_2$ bis $\lambda_h$ darstellen soll.
F"ur $\chi^\lambda(1)$ gilt nach Formel (\ref{Mu-Na}):
\begin{eqnarray*}
|\chi^\lambda(1)| & = & \# \mbox{Abbaum"oglichkeiten, so dass stets } \lambda_i\geq \lambda_j  \ \forall i\geq j \mbox{gilt} \\
& \geq & \# \mbox{Abbaum"oglichkeiten, so dass B am Schluss abgebaut wird} \\
& \geq & \# \mbox{ der M"oglichkeiten A Boxen aus A+B zu w"ahlen} \\
& = & \binom{A+B}{A} \\
& \geq & \binom{A + \frac{1}{2}n}{A} \\
& \geq & \biggl(\frac{{\frac{1}{2}n}}{A}\biggr)^A  \geq 2^A\\
\end{eqnarray*}
Ist nun $A>\sqrt{n}$ so folgt $\chi^\lambda(1) \geq 2^{\sqrt{n}}$ und f"ur $n$ gro"s genug (ca. ab $n$=$8.24*10^8$) gilt $\sqrt{2^{\sqrt{n}}} \geq (2n)^{\frac{\sqrt[3]{n}}{2}}$.

\fbox{{\bf Fall 3: $\lambda_1 \geq \frac{3}{4}n$ und $2\sqrt[3]{n}<A<\sqrt{n}$ }} \\ 
F"ur $A<\sqrt{n}$ ist \\
$\chi^\lambda(1) \geq \biggl(\frac{\frac{1}{2}n}{A}\biggr)^A = \biggl(\frac{\frac{1}{2}n}{A}\biggr)^A$ mit $A>2\sqrt[3]{n}$ folgt dann \\
$\chi^\lambda(1) \geq \biggl(\frac{\sqrt{n}}{2}\biggr)^{2\sqrt[3]{n}}$ und f"ur $n>44$ gilt $\sqrt{2^{\sqrt{n}}} \geq (2n)^{\frac{\sqrt[3]{n}}{2}}$.

\fbox{{\bf Fall 4: $\lambda_1 \geq \frac{3}{4}n$ und $4<A<2\sqrt[3]{n}$ }} \\ 
In diesem Fall l"asst sich $\chi^\lambda(1)^{\frac{1}{2}} \geq \chi^\lambda(\pi')$ nicht durch das Absch"atzen von $\chi^\lambda(1)$ erreichen. Es gilt also $\chi^\lambda(\pi')$ mit $m\leq \sqrt[3]{n}$ besser abzusch"atzen. Sei daf"ur $\pi'=c_1\cdots c_m$, $c_i$ Zykel von $\pi'$. So gilt nach Formel \ref{Mu-Na}
\begin{eqnarray*}
|\chi^\lambda(\pi')| & \leq & \# \mbox{Abbaum"oglichkeiten durch die Zykel} \ c_i\\
& = & \sum_{I\subset \left\{1,\ldots m\right\}} \# \mbox{Abbaum"oglichkeiten mit} \ c_i \ \mbox{ganz in} \ [\mu] \Leftrightarrow i\in I \\
& \leq & \binom{k}{n-\lambda_1} \chi^\mu(1) \ \leq \ \binom{\sqrt[3]{n}}{n-\lambda_1} \chi^\mu(1)
\end{eqnarray*}
F"ur $|\chi^\lambda(1)|$ gilt nach \cite[Lemma 8 (ii)]{pm2007}
\begin{eqnarray*}
|\chi^\lambda(1)| & \geq & \binom{\lambda_1}{n-\lambda_1} \chi^\mu(1) \\
 & = & \binom{2\lambda_1-n+n-\lambda_1}{n-\lambda_1} \chi^\mu(1) \\
 & \geq & \biggl(\frac{2\lambda_1-4}{n-\lambda_1}\biggr)^{n-\lambda_1} \chi^\mu(1) \\
  & \geq & \biggl(\frac{n/2}{2\sqrt[3]{n}}\biggr)^{n-\lambda_1} \chi^\mu(1) \ = \ \biggl(\frac{n^{\frac{2}{3}}}{4}\biggr)^{n-\lambda_1} \chi^\mu(1) 
\end{eqnarray*} 

\textsc{Behauptung:} $\binom{\sqrt[3]{n}}{n-\lambda_1}^2 \chi^\mu(1)^2 \leq \biggl(\frac{n^{\frac{2}{3}}}{4}\biggr)^{n-\lambda_1} \chi^\mu(1) $ \\
Es gilt $\chi^\mu(1) \leq (n-\lambda_1)!^{\frac{1}{2}}$ und $\binom{\sqrt[3]{n}}{n-\lambda_1}^2 \leq \frac{((n^\frac{1}{3})^{(n-\lambda_1)})^2}{(n-\lambda_1)!^2} = \frac{(n^\frac{2}{3})^{(n-\lambda_1)}}{(n-\lambda_1)!^2}$,
das hei"st die Behauptung gilt, da 
\[\begin{array}{cccc}
  & \frac{1}{(n-\lambda_1)!^{\frac{3}{2}}} & \leq & \frac{1}{4^{n-\lambda_1}}\\
\Leftrightarrow &  4^{n-\lambda_1} & \leq & (n-\lambda_1)!^{\frac{3}{2}}\\
 \Leftrightarrow & 4^A & \leq & A!^{\frac{3}{2}} \\
 \Leftrightarrow & 4 & \leq & A \\
\end{array}\]

\fbox{{\bf Fall 5: $\lambda_1 \geq \frac{3}{4}n$ und $A<2\sqrt[3]{n}, A<4$ }} \\ 
Dieser Fall verwendet die Theorie der Charakterpolynome, die in \cite[Kapitel 6.2]{kerber} und in \cite{specht} nachzulesen ist. Sie besagt, dass falls $\pi' \in S_n$ $z_i$ $i-$Zykel besitzt, so gilt f"ur das Charakterpolynom $P$:
\begin{eqnarray*}
P(z_1,\ldots,z_n) & < & c(\sum_{i=1}^n z_i)^{n-\lambda_1}\\
& < & c (\sqrt[3]{n})^{n-\lambda_1}
\end{eqnarray*}
und $P(z_1,\ldots, z_n)=\chi^\lambda(\pi')$.\\
Daraus ergibt sich $\chi^\lambda(\pi')^2 < c^2((\sqrt[3]{n})^2)^{n-\lambda_1}$. F"ur $\chi^\lambda(1)$ gilt in diesem Fall
\begin{eqnarray*}	
\chi^\lambda(1) \geq \binom{\lambda_1}{n-\lambda_1} = \binom{2\lambda_1-n+n-\lambda_1}{n-\lambda_1} \geq \biggl(\frac{2\lambda_1-n}{n-\lambda_1}\biggr)^{n-\lambda_1} \\
= \biggl(\frac{n-2A}{A}\biggr) > \biggl(\frac{n-80}{4}\biggr)^A = \biggl(\frac{n}{4}-2\biggr)^A
\end{eqnarray*}
\newpage
\textsc{Behauptung:} $c^2((\sqrt[3]{n})^2)^{n-\lambda_1} \leq \biggl(\frac{n}{4}-2\biggr)^{n-\lambda_1}$\\
Die Behauptung gilt, da
\[\begin{array}{cccc}
& c'(\sqrt[3]{n})^2 & \leq & \frac{n}{4}-2\\
\Leftrightarrow & c'n^{\frac{2}{3}}+2 & \leq & \frac{n}{4}\\
\Leftrightarrow & 4c'n^{\frac{2}{3}}+8 & \leq & n. \\
\end{array}\]\\ 
Diese letzte Zeile gilt f"ur $n$ abh"angig von $c'$ gro"s genug gew"ahlt. ($c'=1 \Rightarrow n \geq 86, c'=2 \Rightarrow n \geq 536, c'=3 \Rightarrow n \geq 1752$)

Somit gilt in allen f"unf F"allen $\chi^\lambda(1)^{\frac{1}{2}} \geq \chi^\lambda(\pi')$.

\textsc{Teil 2:} \\
Mit der in Teil 1 gezeigten Ungleichung erh"alt man f"ur den Summationsterm $\sum_\lambda \frac{\chi^\lambda
(\pi')}{\chi^\lambda(1)}$ in (\ref{summand}) die folgende Absch"atzung:
\begin{eqnarray*}
\sum_\lambda \frac{\chi^\lambda
(\pi')}{\chi^\lambda(1)} \leq \sum_\lambda \frac{\chi^\lambda (1)^{\frac{1}{2}}}{\chi^\lambda(1)} = \sum_\lambda \frac{1}{\chi^\lambda(1)^{\frac{1}{2}}}
\end{eqnarray*}
Mit Formel (\ref{lulov}) aus Satz 2.2.11 erh"alt man, da die Summation "uber alle irreduziblen Darstellungen $\lambda$ mit $\lambda \neq [n],[1^n],[n-1,1],[2,1^{n-2}]$ l"auft:
\[\sum_\lambda \frac{1}{\chi^\lambda(1)^{\frac{1}{2}}} = \mathcal{O}(n^{-2\cdot \frac{1}{2}}) = \mathcal{O}(n^{-1}).\]
Hieraus folgt nun die Aussage des Theorems, denn nach (\ref{summand}) gilt
\begin{eqnarray*}
P(\sigma\tau=\pi) & = & \frac{1}{n!^4}((1 + sgn(\pi'))(1 + \frac{l-1}{n-1})+\mathcal{O}(n^{-1}))\\
 & = & \frac{1}{n!^4}(1 + sgn(\pi')+ \mathcal{O}(\frac{l}{n^{-1}})) \\
 & = & \frac{1}{n!^4}(1 + sgn(\pi')+ \mathcal{O}(n^{-\frac{2}{3}})).
\end{eqnarray*} \qed 

F"ur die Verteilung der Eckenanzahl ergibt sich nun folgendes:
\begin{eqnarray} \label{summe}
 P(\# \mbox{Ecken} \leq k) & = & \sum_{s \leq k} P(\# \mbox{Ecken} = s) \nonumber \\
& = & \sum_{s \leq k} \sum_{\pi :\ \pi'\ hat \ s \ Zykel \atop s \leq \sqrt[3]{n}} \frac{1}{n!^4}(1 + sgn(\pi')+ \mathcal{O}(n^{-\frac{2}{3}})) + \mathcal{O}(n^{-\frac{1}{3}})
\end{eqnarray}
Der letzte Term ergibt sich leicht durch die Berechnung der Wahrscheinlichkeit $P( \sigma\tau = \pi)$, f"ur $m \geq \sqrt[3]{n}$. Verwendet wird hierf"ur die Gleichung $P(\pi \ \mbox{hat k Zykel } | \pi=\sigma\tau)=(1+\mathcal{O}(n^{-\frac{1}{3}})P(\pi \ \mbox{hat k Zykel} | U),$
wobei $U$ die Gleichverteilung ist.
\section{Verteilung des Geschlechts}
{\bf Theorem 3.3.1} \begin{it}\\
Sei $X$ die Zufallsvariable, die jedem zuf"allig gew"ahlten Origami das Geschlecht $g$ zuordnet, $E_n$ der Erwartungswert und $\sigma_n$ die Standardabweichung der Zufallsvariablen $X$, so erh"alt man f"ur $\sigma_n$ und $E_n$ die folgenden Werte
\[E_n=-\frac{1}{2}(\log n + \gamma + o(1))+1+\frac{1}{2}n\]
\[\sigma_n = \frac{1}{2}(\sqrt{\log n} - (\frac{\pi^2}{12}-\frac{\gamma}{2})\frac {1}{\sqrt {\log n}} + o(\frac {1}{\sqrt {\log n}}))
\] wobei $\gamma$ die Eulerkonstante ist und die Zufallsvariable $X$ ist normalverteilt.
\end{it} \\
\bew
Nach Gleichung (\ref{summe}) gilt
\begin{eqnarray*} 
P(\# \mbox{Ecken} \leq k) & = & \frac{1}{n!^4} \sum_{s \leq k} \sum_{\pi :\ \pi'\ hat \ s \ Zykel \atop s \leq \sqrt[3]{n}} (1 + sgn(\pi')+ \mathcal{O}(n^{-\frac{2}{3}})) + \mathcal{O}(n^{-\frac{1}{3}}).
\end{eqnarray*}
Es wird jetzt der zweite Summand nach dem Kriterium $\pi \in A_n$ oder $\pi \in S_n\setminus A_n$ gesplittet, um danach diese beiden Teilsummen genauer zu betrachten.
\begin{eqnarray*}
\sum_{\pi :\ \pi'\ hat \ s \ Zykel \atop s \leq \sqrt[3]{n}} 1 + sgn(\pi')+\mathcal{O}(n^{-\frac{2}{3}})= \underbrace{\sum_{\pi :\ \pi'\ hat \ s \ Zykel \atop s \leq \sqrt[3]{n}, \pi' \in S_n}}_A 2 + \mathcal{O}(n^{-\frac{2}{3}})+\underbrace{\sum_{\pi :\ \pi'\ hat \ s \ Zykel \atop s \leq \sqrt[3]{n}, \pi' \in S_n\setminus A_n}}_B \mathcal{O}(n^{-\frac{2}{3}})
\end{eqnarray*}
Die Summen $A$ und $B$ lassen sich wie folgt darstellen, dabei sei daran erinnert, dass $\pi'=\pi_1\circ\pi_2\circ\pi_3\circ\pi_4$, $\pi_i \in S_n \ \forall i$ ist.
\begin{eqnarray*}
A & = & \biggl(\frac{1}{2} n! n! n! \biggl[ \# \underbrace{\biggl\{\pi_4 \in A_n | \pi'\ hat \ s \ Zykel \biggr\}}_\alpha + \# \underbrace{\biggr\{\pi_4 \in S_n\setminus A_n | \pi'\ hat \ s \ Zykel \biggr\}}_\beta\biggr]\biggr) \\
B & = & \biggl(\frac{1}{2} n! n! n! \biggl[ \# \alpha + \# \beta \biggr]\biggr)
\end{eqnarray*}
Daraus folgt
\begin{eqnarray*}
P(\# \mbox{Ecken} \leq k) & = &  \frac{1}{n!^4}(2+\mathcal{O}(n^{-\frac{2}{3}})) \sum_{s \leq k} A + \frac{1}{n!^4} \mathcal{O}(n^{-\frac{2}{3}}) \sum_{s \leq k} B + \mathcal{O}(n^{-\frac{1}{3}}) \\
& = & \frac{1}{2} \biggl(2+\mathcal{O}(n^{-\frac{2}{3}})\biggr) \biggl( \sum_{s \leq k} \frac{1}{n!} \biggl[ \# \alpha + \# \beta\biggr]\biggl) \\
&  & + \frac{1}{2}\mathcal{O}(n^{-\frac{2}{3}}) \biggl( \sum_{s \leq k} \frac{1}{n!} \biggl[ \# \alpha + \# \beta \biggr]\biggl)+ \mathcal{O}(n^{-\frac{1}{3}})
\end{eqnarray*}
Nach Theorem 2.2.12 gilt, dass die Anzahl $M_{(n)}$ der Zykel einer Permuation normalverteilt ist. Die Normalverteilung wird hier mit $\mathcal{N}(E_n,\sigma_n^2)$ bezeichnet. Somit ergibt sich f"ur die Verteilung der Eckenanzahl
\begin{eqnarray*}
P(\# \mbox{Ecken} \leq k) & = & \biggl(1 + \frac{1}{2}\mathcal{O}(n^{-\frac{2}{3}})\biggr) \biggl[ \frac{1}{2}\mathcal{N}(E_n,\sigma_n^2) + \frac{1}{2}\mathcal{N}(E_n,\sigma_n^2)\biggr] \\
& & + \frac{1}{2}\mathcal{O}(n^{-\frac{2}{3}})\biggl[ \frac{1}{2}\mathcal{N}(E_n,\sigma_n^2) + \frac{1}{2}\mathcal{N}(E_n,\sigma_n^2)\biggr] + \mathcal{O}(n^{-\frac{1}{3}}) \\
& = & \mathcal{N}(E_n,\sigma_n^2) + \mathcal{O}(n^{-\frac{2}{3}}) \mathcal{N}(E_n,\sigma_n^2) + \mathcal{O}(n^{-\frac{1}{3}})
\end{eqnarray*}
Sei $Y$ die Zufallsvariable, die jeder Verklebungsvorschrift die daraus resultierende Anzahl der Ecken zuordnet, so besteht zwischen der Zufallsvariablen $Y$ und der Zufallsvariablen $X$ die Beziehung:
\[ X= \frac{1}{2} (-Y+2+K-F) = -\frac{1}{2}Y+1+\frac{1}{2}n\]
Das hei"st $X$ ergibt sich aus $Y$ durch lineare Transformation. F"ur den Erwartungswert von $X$ gilt dann 
\[ E(X)= -\frac{1}{2}E(Y)+1+\frac{1}{2}n\]
und f"ur die Standardabweichung gilt
\[ \sigma_X=Var(X)^{\frac{1}{2}}=((-\frac{1}{2})^2 Var(Y))^\frac{1}{2}=\frac{1}{2} \sigma_Y.\]
Mit Theorem 2.2.12 folgt somit die Aussage von Theorem 3.3.1. \qed
\section{Charakterberechnung f"ur Elemente aus $H\wr S_n$}
Das wesentliche Resultat dieses Unterkapitels wird sein, die Charakterwerte f"ur Elemente $(h;f)$ aus $H\wr S_n$ mit der Eigenschaft, dass die Elemente $f_i \in S_n$ $r$-Zykel sind, abzusch"atzen. Diese Absch"atzung wird nach (\ref{3}) zur"uckgef"uhrt auf die Absch"atzung von $\chi^\lambda$, wobei $[\lambda]$ irreduzible Darstellung von $S_n$ ist. Hierf"ur wird zun"achst das f"ur die Absch"atzung der Charakterwerte $\chi^\lambda$ wesentliche Theorem von Fomin und Lulov angegeben.

{\bf Theorem 3.4.1} \\ \begin{it}
Sei $r$ eine fixe ganze Zahl $\geq 2$, $n$ durch $r$ teilbar, n=mr, und sei $\pi \in S_n$ eine Permutation vom Zykeltyp $(r^m)$. Dann gilt f"ur alle irreduziblen Darstellungen $[\lambda]$
\begin{eqnarray}\label{fl}
|\chi^\lambda(\pi)| \leq \frac{m! \ r^m}{(n!)^{\frac{1}{r}}} \cdot \chi^\lambda(1)^{\frac{1}{r}} \leq c \cdot n^{\frac{1}{2}(1-\frac{1}{r})} \cdot \chi^\lambda(1)^{\frac{1}{r}},
\end{eqnarray}
wobei c nur von r abh"angt.\end{it}\\
\bew siehe \cite{fm}. \qed

Wie man sp"ater im Beweis von Theorem 3.4.3 sehen wird, l"asst sich obiges Theorem im allgemeinen Fall, das hei"st $h \in H$ beliebig gew"ahlt, nicht anwenden. Man ben"otigt dazu eine Absch"atzung f"ur Elemente des Zykeltyps $(1^b,r^m)$, die sich jedoch aus Theorem 3.4.1 und der Murnaghan-Nakayama Formel gewinnen l"asst.

{\bf Lemma 3.4.2} \\ \begin{it}
Sei $\pi \in S_n$ vom Zykeltyp $(1^b,r^m)$, dann gilt f"ur alle irreduziblen Darstellungen $\lambda$ von $S_n$
\begin{eqnarray}\label{fl2}
|\chi^\lambda(\pi)| \leq c \cdot (2n)^{\frac{1}{2}(b+1)} \cdot \chi^\lambda(1)^{\frac{1}{r}},
\end{eqnarray}
wobei c nur von r abh"angt.\end{it}\\
\bew siehe \cite[S.568,Prop 2.12]{ls2004}.\qed 

{\bf Theorem 3.4.3} \\ \begin{it}
Sei $H$ eine Permutationsgruppe auf der Menge $\Omega$ mit $|\Omega|=k$, $S_n$ die symmetrische Gruppe und sei $(h;f) \in H\wr S_n$ von der Form $(h;f_1,\ldots ,f_k)$, wobei die Elemente $f_i \in S_n$ $r$-Zykel sind und $h$ vom Zykeltyp $Th=(a_1,\ldots ,a_k)$ ist. Dann gilt f"ur alle $D=D'\otimes \widetilde{\lambda^*} \in$ Irr($H\wr S_n$)
\[|\chi^D(h;f)| \leq \chi^{D^{''}}(h) \cdot c^{\sum_{i=1}^k a_i} \cdot (2n)^{\frac{1}{2}(\sum_{i=1}^k b_i a_i +a)} \cdot \prod_{i=0}^{k-1} \prod_{\nu=1}^{a_{i+1}} \chi^{\lambda_{j_{\sum_{l=0}^i a_l+\nu}}}(1)^{\frac{1}{r}},\] wobei $c$ nur von $r$ abh"angt. \end{it} \\
\bew Nach (\ref{3}) gilt 
\begin{eqnarray*}
\chi^{D'\otimes \widetilde{D^*}}(h;f)&=&\chi^{D''}(h)\prod_{\nu=1}^{\sum a_i} \chi^{\lambda_{j_\nu}}(g_\nu(h;f)) \\
&=& \chi^{D''}(h) \ \chi^{\lambda_{j_1}}(\sigma_{1,1}) \cdots \chi^{\lambda_{j_{a_1}}}(\sigma_{1,a_1}) \cdot \chi^{\lambda_{j_{a_1+1}}}(\sigma_{2,1})   \cdots \chi^{\lambda_{j_{a_1+a_2}}}(\sigma_{2,a_2}) \cdots \\ & & \cdot \chi^{\lambda_{j_{\sum_{i=1}^{k-1} a_i}}}(\sigma_{k,1}) \cdots \chi^{\lambda_{j_{\sum a_i} }}(\sigma_{k,a_k}) \\
&=& \chi^{D''}(h) \prod_{i=0}^{k-1} \prod_{\nu=1}^{a_{i+1}} \chi^{\lambda_{j_{\sum_{l=0}^i a_l+\nu}}}(\sigma_{i,\nu}),
\end{eqnarray*}
wobei die Elemente $\sigma_{i,\nu} \in S_n$ vom Zykeltyp $(1^{b_i},r^i)$ sind und $a_0:=0$ ist.\\
Wendet man Lemma 3.4.2 auf die Charaktere $\chi^{\lambda_{j_s}}$, $s=\sum_{l=0}^i a_l+\nu$, an, dann folgt
\begin{eqnarray*}
\chi^{D'\otimes \widetilde{D^*}}(h;f) & \leq & \chi^{D''}(h) \prod_{i=0}^{k-1} \prod_{\nu=1}^{a_{i+1}} \cdot c \cdot (2n)^{\frac{1}{2}(b_i+1)} \chi^{\lambda_{j_s}}(1)^{\frac{1}{r}}  \\
& = & \chi^{D''}(h) \prod_{i=0}^{k-1} \cdot c^{a_{i+1}} \cdot (2n)^{\frac{1}{2}(b_i+1)^{a_{i+1}}} \prod_{\nu=1}^{a_{i+1}} \chi^{\lambda_{j_s}}(1)^{\frac{1}{r}}  \\
& = & \chi^{D''}(h) \ c^{\sum a_i} \cdot (2n)^{\frac{1}{2}(\sum_{i=1}^k b_i a_i +\sum a_i)} \prod_{i=0}^{k-1} \prod_{\nu=1}^{a_{i+1}} \chi^{\lambda_{j_s}}(1)^{\frac{1}{r}}. \hspace{2.5cm} \qed
\end{eqnarray*} 

{\bf Korollar 3.4.4} \\ \begin{it}
Unter den Voraussetzungen des Theorem 3.4.3 gilt
\[\chi^{D'\otimes \widetilde{D^*}}(h;f) \leq c^k \cdot (2n)^{\frac{k}{2}(n-r+1)} \cdot \chi^{D^{''}}(1)^{\frac{r-1}{r}} \cdot \chi^D(1)^{\frac{1}{r}}.\]
\end{it}
\bew Es gilt $\sum_{i=1}^k a_i \leq k, \sum_{i=1}^k i \cdot a_i=k$ und $b_i=n-i \cdot r$, es ergibt sich somit die Absch"atzung
\[\frac{1}{2}(\sum_{i=1}^k b_i a_i +a) \leq \frac{k}{2}(n-r+1).\]
Nach Gleichung (\ref{3}) gilt f"ur $\chi^D(1)$, 
\begin{eqnarray*}
\chi^D(1_H;e)= \chi^{D^{''}}(1) \prod_{\nu=1}^k \chi^{\lambda_{j_\nu}}(1), 
\end{eqnarray*} 
woraus folgt, dass 
\begin{eqnarray*}
\prod_{i=0}^{k-1} \prod_{\nu=1}^{a_{i+1}} \chi^{\lambda_{j_{\sum_{l=0}^i a_l+\nu}}}(1)^{\frac{1}{r}} &\leq&  \prod_{\nu=1}^k \chi^{\lambda_{j_\nu}}(1)^{\frac{1}{r}} \\ &=& \frac{\chi^D(1)^{\frac{1}{r}}}{\chi^{D^{''}}(1)^\frac{1}{r}}.
\end{eqnarray*} 
Es gilt
\begin{eqnarray*}
\chi^{D'\otimes \widetilde{D^*}}(h;f) &\leq& \chi^{D''}(h) \cdot c^k \cdot (2n)^{\frac{k}{2}(n-r+1)} \cdot \frac{\chi^D(1)^{\frac{1}{r}}}{\chi^{D^{''}}(1)^\frac{1}{r}} \\
&\leq& c^k \cdot (2n)^{\frac{k}{2}(n-r+1)} \cdot \chi^{D^{''}}(1)^{\frac{r-1}{r}} \cdot \chi^D(1)^{\frac{1}{r}}.
\end{eqnarray*} \qed

\chapter{Origamis in der Forschung}
Dieses Kapitel soll eine kurze Einf"uhrung der Begriffe Translationsstruktur, Teichm"uller Kurve im Modulraum und Veech Gruppe geben, da sich die Origami-Forschung vorwiegend mit Origamis als Teichm"uller Kurven besch"aftigt. Daf"ur werden auch die in 1.1 angef"uhrten Beschreibungen von Origamis explizit angegeben. Die Einf"uhrung geht vor allem auf die Arbeiten der Forschungsgruppe der Universtit"at Karlsruhe zur"uck. Diese besch"aftigt sich schon seit einigen Jahren mit Origamis und arbeitet seit Sommer 2007 unter dem Forschungsprojekt "'Mit Origamis zu Teichm"ullerkurven im Modulraum"'. Andere Autoren untersuchen Origamis unter dem Namen "'square tiled surfaces"', die zur allgemeineren Gruppe der "'flat surfaces"' geh"oren. Diese Gruppe wird schon seit l"angerer Zeit in den Teilgebieten der algebraischen Geometrie, der komplexen Analysis und der dynamischen Systeme untersucht.

Die in den n"achsten beiden Kapiteln verwendeten Grundlagen der "Uberlagerungstheorie k"onnen zum Beispiel in \cite{st} nachgelesen werden.
\section{Origamis als "Uberlagerungen des Torus}
Sei $X$ die zusammenh"angende und abgeschlossene Fl"ache aus Definition 1.1.1 und $T$ der durch das Origami $O_0$ enstehende Torus, so existiert eine nat"urliche Abbildung $X \rightarrow T$, indem man jedes Einheitsquadrat des Origami $O$ auf das Einheitsquadrat des Origami $O_0$ abbildet. Diese Abbildung ist eine "Uberlagerung, die au"ser "uber dem Eckpunkt $P$ unverzweigt ist. Umgekehrt erh"alt man zu jeder Abbildung $p: X \rightarrow T$ einer abgeschlossenen Fl"ache $X$, eine Zerlegung von $X$ in Quadrate mittels Zerschneidung von $X$ entlang der Urbilder der Ecken des Origami $O_0$.

{\bf Definition 4.1.1.1} \\
Ein {\it Origami $O$ vom Geschlecht $g\geq1$} und {\it vom Grad $n$} ist eine "Uberlagerung $p: X \rightarrow T$ einer abgeschlossenen, orientierten Fl"ache $X$ vom Geschlecht $g$ in den Torus $T$, die nur "uber einem Punkt $P \in T$ verzweigt. 

{\bf Definition 4.1.1.2} \\
Zwei Origamis $O_1 = (p_1: X_1 \rightarrow T)$ und $O_2 = (p_2: X_2 \rightarrow T)$ hei"sen
{\it "aquivalent}, falls ein Homeomorphismus $\varphi: X_1 \rightarrow X_2$ existiert, so dass $p_1=p_2 \circ \varphi$ gilt.
\section{Origamis als Untergruppe der freien Gruppe $\mathbb{F}_2$ mit endlichem Index}
Sei $O = (p:X \rightarrow T)$ ein Origami. Sei $T^*=T - P$ und $X^*=X - p^{-1}(P)$, dann ist $p:X^* \rightarrow T^*$ eine unverzweigte "Uberlagerung. Dies f"uhrt zu einer Einbettung der korrespondierenden Fundamentalgruppen:
\[ U=\pi_1(x^*) \subseteq \pi_1(T^*) \cong \mathbb{F}_2.\]
Dabei erh"alt man den fixierten Isomorphismus $\pi_1(T^*) \cong \mathbb{F}_2$, indem man $M$ als Mittelpunkt von $T$ und die einfache abgeschlossene horizontale bzw. vertikale Kurve durch $M$ als Erzeuger der Fundamentalgruppe $\pi_1(T^*,M)$ festlegt. $\mathbb{F}_2$ ist dabei die freie Gruppe mit den zwei Erzeugern $x$ und $y$, siehe Abbildung 5.
\begin{center}
\begin{picture}(2,0)
\put(0,0){\line(1,0){2}}
\put(0,0){\line(0,-1){2}}
\put(2,0){\line(0,-1){2}}
\put(0,-2){\line(1,0){2}}
\put(0,-1){\vector(1,0){2}}
\put(1,-2){\vector(0,1){2}}
\put(0.9,-1.1){$\bullet$}
\put(1.2,-0.9){$M$}
\put(0.5,-0.9){$x$}
\put(1.2,-1.6){$y$}
\end{picture}
\end{center}
\vspace{1cm}
\begin{center}
Abb.5: Erzeuger von $\pi_1(T^*)$.
\end{center}
W"ahlt man ein zu $O$ "aquivalentes Origami, so f"uhrt dies zu einer Konjugation von $U$ mit einem Elemente aus $\mathbb{F}_2$. Der Index der Untergruppe von $\mathbb{F}_2$ ist der Grad $n$ der "Uberlagerung $p$. Ist andererseits eine Untergruppe $U$ von endlichem Index von $\mathbb{F}_2$ gegeben, so erh"alt man wie folgt mit dem Theorem der universellen "Uberlagerung eine "Uberlagerung $X \rightarrow T$.\\
Sei $q: \tilde{T}^* \rightarrow T^*$ die universelle "Uberlagerung von $T^*$, so ist nach dem Theorem der universellen "Uberlagerung $\pi_1(T^*)$ isomorph zu Deck($\tilde{T}^*/T^*$), der Gruppe der Decktransformationen von $\tilde{T}^*/T^*$, und jede Untergruppe $U \subseteq$ Deck($\tilde{T}^*/T^*$) von endlichem Index  korrespondiert mit einer unverzweigten "Uberlagerung $p: \tilde{X}^* \rightarrow T^*$ von endlichem Grad. Diese kann zu einer "Uberlagerung $X \rightarrow T$, wobei X eine abgeschlossene Fl"ache ist, erweitert werden. Daraus folgt die 1-1 Korrespondenz: 

\hspace{1cm} Origamis \hspace{2.5cm} $\leftrightarrow$ \hfill Untergruppe von $\mathbb{F}_2$ von endlichem Index \\
\hspace{0.5cm}bis auf Isomorphismus \hspace{5.5cm} bis auf Konjugation.
\section[Translationsstruktur, Teichm"uller Kurve im Modulraum \\ \hspace*{0.1cm}und Veech Gruppe]{Translationsstruktur, Teichm"uller Kurve im \\ Modulraum und Veech Gruppe}
\subsection{Translationsstruktur}
Ein Origami $O = (p:X \rightarrow T)$ definiert auf nat"urliche Weise eine $SL_2(\mathbb{R})$-Familie von Translationsstrukturen $\mu_A$, $A \in SL_2(\mathbb{R})$, auf $X^*$. Dabei hei"st eine Struktur, d.h. eine durch einen Atlas bestimmte "Aquivalenzklasse, {\it Translationsstruktur}, falls alle Karten"uberg"ange des Atlas Translationen sind.
\\
Sei $A=\begin{pmatrix} a & b \\ c & d \end{pmatrix} \in SL_2(\mathbb{R})$, so definiert diese Matrix eine Translationsstruktur $\eta_A$ auf dem Torus $T$, indem dieser mit $\mathbb{C}/\Lambda_A$ identifiziert wird; dabei ist $\Lambda_A$ das Gitter $\left\langle \binom{a}{c}, \binom{b}{d} \right\rangle$ in $\mathbb{C}$. Die Translationsstruktur $\mu_A$ auf $X^*$ ist dann die Hochhebung der Translationsstrukur $\eta_A$ durch $p$, d.h. $\mu_A=p^*\eta_A$.

Die anf"angliche Beschreibung von Origamis durch Verklebung der Kanten ergibt die Translationsstruktur $\mu_I$, $I$ Einheitsmatrix, indem man die Quadrate mit dem euklidischen Einheitsquadrat in $\mathbb{C}$ identifiziert. F"ur eine beliebige Matrix $A \in SL_2(\mathbb{R})$ erh"alt man $\mu_A$ durch die Identifikation der Quadrate mit dem Parallelogramm, das durch die beiden Vektoren $\binom{a}{c}$ und $\binom{b}{d}$ aufgespannt wird.
\subsection{Teichm"uller Kurve im Modulraum}
Der {\it Modulraum} $M_g$ einer abgeschlossenen Riemannschen Fl"ache $R$ vom Geschlecht $g$ ist der Quotient $T_g/Mod(R)$. $T_g$ ist der Teichm"uller Raum von $R$ und $Mod(R)$ die Menge der Homotopieklassen $[\omega]$, der orientierungserhaltenden Diffeomorphismen $\omega: R \rightarrow R$.\\
Die hier verwendeten Begriffe, speziell der des Teichm"uller Raums und der unten verwendete Begriff der Markierung, k"onnen zum Beispiel in \cite[S.14]{ima} nachgelesen werden.
 
Ein Origami $O = (p:X \rightarrow T)$ definiert eine {\it Teichm"uller Kurve im Modulraum} $M_g$ einer abgeschlossenen Riemannschen Fl"ache vom Geschlecht $g$.\\
Betrachte dazu die oben konstruierte Translationsstruktur $\mu_A$ auf $X^*$. Diese kann auf die abgeschlossene Fl"ache $X$ erweitert werden. Die Riemannsche Fl"ache $(X,\mu_A)$ zusammen mit der Identit"atsabbildung als Markierung definiert dann einen Punkt im Teichm"uller Raum $T_g$, so dass man die folgende Abbildung erh"alt:
\[ \hat{\iota}: SL_2(\mathbb{R}) \rightarrow T_g, \ A \mapsto [(X,\mu_A),id].\]
Falls $A \in SO_2(\mathbb{R})$, so ist die affine Abbildung $z \mapsto A \cdot z$ holomorph und die Abbildung $\hat{\iota}$ faktorisiert "uber $SO_2(\mathbb{R})$. Verwendet man weiter, dass $SL_2(\mathbb{R})/SO_2(\mathbb{R})\cong \mathbb{H}$, so erh"alt man die Abbildung 
\[ \iota: \mathbb{H} \rightarrow T_g.\]
Das Bild $\Delta$ dieser holomorphen und isometrischen {\it Teichm"uller Einbettung} hei"st {\it Teich\-m"ullerkreisscheibe} oder {\it Geod"atische}.
Betrachtet man das Bild solch einer Teichm"ullerkreisscheibe unter der Projektion $\pi$ von $T_g$ auf den Modulraum $M_g$, so ist dies eine komplexe algebraische Kurve und alle Kurven in $M_g$, die als Bilder von Teichm"ullerkreisscheiben auftreten, nennt man {\it Teichm"uller Kurven}. Betrachtet man die zu einem Origami $O$ geh"orige komplexe Geod"atische $\Delta$, so hei"st $\pi(\Delta)$ {\it Origamikurve}.
\subsection{Veech Gruppe}
Sei $X^*$ zusammenh"angede Fl"ache und $\mu$ Translationsstruktur auf $X^*$, so l"asst sich eine dazugeh"orige Untergruppe von $GL_2(\mathbb{R})$ bestimmen. Diese nennt man {\it Veech Gruppe}.

Dazu sei\[Aff^+(X^*,\mu):=\left\{f: X^* \rightarrow X^* | \mbox{ f orientierungserhaltender affiner Diffeomorphismus}\right\},\] wobei der affine Diffeomorphismus $f \in  Aff^+(X^*,\mu)$ lokal gegeben ist durch die Abbildung $z \mapsto A \cdot z + b$ mit $A \in SL_2(\mathbb{R})$ und $b \in  \mathbb{R}^2$. $A$ h"angt dabei nur von $f$ ab und ist unabh"angig von der Wahl der Karten, so dass man die wohldefinierte Abbildung
\[ D: Aff^+(X^*,\mu) \rightarrow SL_2(\mathbb{R}), \ f \mapsto A \]
erh"alt. Diese nennt man {\it Derivationsabbildung}. Die Veech Gruppe $\Gamma(X^*,\mu)$ der Translationsfl"ache $(X^*,\mu)$ ist dann das Bild der Derivationsabbildung D: 
\[ \Gamma(X^*,\mu) = D(Aff^+(X^*,\mu))\]
Sei nun $O = (p:X \rightarrow T)$ ein Origami und $\mu_A$ die durch $O$ definierte $SL_2(\mathbb{R})$-Familie von Translationsstrukturen auf $X^*$, so sind die dazugeh"origen Veech Gruppen konjugiert, d.h. $\Gamma(X^*,\mu_A)= A \Gamma(X^*,\mu_I) A^{-1}$. Die Veech Gruppe $\Gamma(O)$ eines Origamis $O$ ist somit definiert durch $\Gamma(X^*,\mu_I)$.

Zum Schluss sei erw"ahnt, dass Gutkin und Judge in \cite{gj} zeigten, dass eine Veech Gruppe genau dann Untergruppe von endlichem Index von $SL_2(\mathbb{Z})$ ist, wenn sie Veech Gruppe eines Origamis ist.

\end{document}